\documentclass[11pt]{article}
\usepackage{CJK}
\usepackage{ws-rotating}
\usepackage{epstopdf}
\usepackage{epsfig}
\usepackage{color}
\usepackage{amsmath}
\usepackage{amssymb}
\usepackage{paralist}
\usepackage{graphicx}
\usepackage{epstopdf}
\usepackage[numbers,sort&compress]{natbib}

\newtheorem{mythm}{\bf{Theorem}\rm}
\newtheorem{myprop}{\bf{Proposition}\rm}
\newtheorem{mylem}{\bf{Lemma}\rm}

\newtheorem{mydef}{\bf{Definition}\rm}

\begin{document}
\pagenumbering{arabic}

\title{Bifurcation of the traveling wave solutions in a perturbed $(1 + 1)$-dimensional dispersive long wave equation via a geometric approach
}
\author{Hang Zheng$^{1,2}$\,\,\,\, Yonghui Xia$^{1}$\footnote{Corresponding author:
 Email: xiadoc@163.com; yhxia@zjnu.cn.}
  \\
{\small \em 1. Department of Mathematics, Zhejiang Normal University, Jinhua, 321004, China}\\
{\small   \,\,\,  yhxia@zjnu.cn; xiadoc@163.com}\\
{\small \em 2.  Department of Mathematics and Computer, Wuyi University, Wuyishan,  354300, China}\\
{\small  \,\,\,  zhenghang513@zjnu.edu.cn; zhenghwyxy@163.com}\\
}

\date{}
\maketitle \pagestyle{myheadings} \markboth{}{}

\noindent

\begin{abstract}
{ Choosing ${\kappa}$ (horizontal ordinate of the saddle point associated to the homoclinic orbit) as bifurcation parameter,  bifurcations of the travelling wave solutions is studied in a perturbed $(1 + 1)$-dimensional dispersive long wave equation. %The solitary wave solution and kink (anti-kink) wave solutions exist in the different bifurcation parameter region, respectively.
The solitary wave solution exists at a suitable wave speed $c$ for the bifurcation parameter ${\kappa}\in (0,1-\frac{\sqrt3}{3})\cup (1+\frac{\sqrt3}{3},2)$, while the kink and anti-kink wave solutions exist at a unique wave speed $c^*=\sqrt{15}/3$ for $\kappa=0$ or $\kappa=2$.%Moreover, the solitary and kink (anti-kink) wave solutions can't coexist at the same bifurcation parameter value.
}
The methods are based on the geometric singular perturbation (GSP, for short) approach, Melnikov method and invariant manifolds theory. {
 Interestingly,
not only the explicit analytical expression of the complicated homoclinic Melnikov integral is directly obtained for the perturbed long wave equation, but also the explicit analytical expression of the limit wave speed is directly given.   Numerical simulations are utilized to verify our mathematical results.}
%Persistence of the solitary wave solution is studied for a perturbed $(1 + 1)$-dimensional dispersive long wave equation. The methods are based on the geometric singular perturbation (GSP, for short) approach, Melnikov method, invariant manifolds and bifurcation analysis. The results show that, under the small singular perturbation, the solitary wave solution does persist at a suitable wave speed $c$ and bifurcation parameter ${\kappa}$ (horizontal ordinate of the saddle point associated to the homoclinic orbit).

%Most of the literature emphasized the existence of the solitary wave by the GSP approach, but there is few literature applying GSP approach to discuss the nonexistence of kink (anti-kink) wave solutions. We prove the nonexistence of kink (anti-kink) wave solutions by GSP approach.
\end{abstract}



\section{Introduction}

\subsection{Model formulation}
%The study of water wave equations (eg. see KdV equation, BB) equation, D) equation and CH or rotation-CH equation)
%have attracted a lot of interest among scholars (see, e.g.,  %\cite{Constantin-acta,Constantin-JFA,fy-JFA,fl-AM,Li-Book,GuiGL,Gui2,Gui3,Luo,CHUjf-JDE2020,CHUjf-SPAM,Ma-Feng-JPA,EJ-JGM,EJ-IUMJ}).
%For instance, Korteweg-de Vries (KdV) equation, Benjamin-Bona-Mahony (BBM) equation, Degasperis-Procesi (DP) equation and Camassa-Holm (CH) equation, etc.

The water wave equations (e.g., KdV equation, BBM equation and CH or rotation-CH equation)
have attracted a lot of scholars' attentions.
These nonlinear wave equations have been extensively used to describe dynamical behaviour of nonlinear waves in shallow water.  Since the interaction of nonlinear and dispersion factors, long wave in shallow water would admit many special characteristics. One of the important property of long waves, is that, they retain their shapes and forms after mutual interactions and collisions. In 1996, Wu and Zhang \cite{wuzhang-book} derived an equation describing nonlinear dispersive long gravity waves travelling in two horizontal directions on shallow waters of uniform depth format formulated as
\begin{equation}\label{wuzhang equation}
\begin{cases}
\theta_t+\theta\theta_x+u\theta_y+v_x=0,\\
u_t+\theta u_x+uu_y+v_y=0,\\
v_t+(\theta v)_x+(uv)_y+\frac{1}{3}(\theta_{xxx}+\theta_{xyy}+u_{xxy}+u_{yyy})=0.
\end{cases}
\end{equation}
where $\theta$ (resp., $u$) is the surface velocity of water along the $x$ (resp., $y$) direction and $v$ is the elevation of the water wave.
For dispersive long wave equation \eqref{wuzhang equation}, Chen, Tang and Lou \cite{FE-PLA2} obtained a special type of multisoliton solution by the Weiss-Tabor-Carnvale
Painlev\'e truncation expansion.
 Eq. (\ref{wuzhang equation}) reduces to the following $(1+1)$-dimensional dispersive long wave equation (for short, DLWE) by symmetry reduction and scale transition:
\begin{equation}\label{wuzhang system}
\begin{cases}
u_t=-uu_x-v_x,\\
v_t=-vu_x-uv_x-\frac{1}{3}u_{xxx}.
\end{cases}
\end{equation}
Due to it models the nonlinear water wave availably, $(1+1)$-dimensional DLWE (\ref{wuzhang system}) is often applied to coastal design and harbor construction. So far, a lot of articles have been concerned with {the exact solutions of $(1+1)$-dimensional DLWE (\ref{wuzhang system})} by using various methods, such as
 a modified Conte's invariant Painlev\'e expansion (Chen and Lou \cite{CL-CSF}), a new Jacobi elliptic function rational expansion method
 (Wang, Chen and Zhang \cite{wq-CSF}), a new general algebraic method with symbolic computation (Chen and Wang \cite{ChenY-AMC}), a generalized extended rational expansion method (Zeng and Wang \cite{zx-AMC}) and extended tanh-function method (Fan \cite{FE-PLA2}).

However, in the real world application, particularly for nonlinear wave, there are many factors of uncertainty and unpredictable perturbations. Thus, small perturbation terms are usually added to describe such unpredictable perturbations when modelling the nonlinear waves. For shallow water wave equation, weak backward diffusion $u_{xx}$ and dissipation $u_{xxxx}$ are called Kuramoto-Sivashinsky (KS, for short) perturbation.
However, to the best of our knowledge, there are no papers considering the perturbed $(1+1)$-dimensional DLWE. In this paper, we study a perturbed $(1+1)$-dimensional DLWE described by
\begin{equation}\label{wuzhang1}
\begin{cases}
u_t=-uu_x-v_x,\\
v_t=-vu_x-uv_x-\frac{1}{3}u_{xxx}+\varepsilon (u_{xx}+u_{xxxx}),
\end{cases}
\end{equation}
where $0<\varepsilon\ll 1$ (a sufficiently small parameter), $u_{xx}$ is backward diffusion and $u_{xxxx}$ represents a dissipation term.

\subsection{Methods, motivation and contributions}

\indent
Note that the perturbed $(1+1)$-dimensional DLWE \eqref{wuzhang1} is singular perturbation system due to the introduction of $\varepsilon u_{xxxx}$.
Considering the geometric singular perturbations in the evolution equation is a effective method to describe the real situation for the nonlinear waves. One of important themes of shallow water waves is to study the existence of traveling wave solutions.
In this paper, we will apply the geometric singular perturbation theory to study the traveling wave solutions of the singular perturbed DLWE \eqref{wuzhang1}.

In fact, the theory of the geometric singular perturbation have been well developed by (Fenichel \cite{FN-JDE}, Jones \cite{JCK-berlin}, Szmolyan \cite{PS-JDE}, Guckenheimer and Holmes \cite{GJ-NY}, Wiggins \cite{WS-NY}, Zhang \cite{ZX-SCM} and Li et al. \cite{LJ-TAMS,LJ-DCDS,LJ-JDE}). And it is a powerful tool to solve the applications arising in nonlinear waves, biological systems and other dynamic systems.
Up till now, there are an extensive literature applying GSP theory  to study such systems, including delayed CH equation (Du et al. \cite{DLL2018}), CH Kuramoto-Sivashinsky equation (Du et al. \cite{DZ-JDE}), generalized CH equation (Qiu et al. \cite{QHM-CNSNS}), perturbed BBM equation (Chen et al. \cite{CAY-JDE}), generalized BBM equation (Sun et al. \cite{SY2019}, Zhu et al. \cite{ZK-ND}),
the delayed DP equation (Cheng and Li \cite{cheng-dcds}),
KdV equation (Derks and Gils \cite{DG-JJIAM}, Ogama \cite{OT-HMJ}), generalized KdV equation (Chen et al. \cite{CAY-AML}, Yan et al. \cite{YW-MMA}, Zhang et al. \cite{Zhanglijun-ND}), biological model (Chen and Zhang \cite{Zhangxiang-SpM}, Wang and Zhang \cite{Zhangxiang-JDE}), reaction-diffusion equation (Yang and Ni \cite{Nimingkang-SCM}, Wu and Ni \cite{Nimingkang-CNSNS}), Belousov-Zhabotinskii system (Du and Qiao \cite{DQ2020}), MEMS model (et al. Iuorio et al \cite{PS2}), piecewise-smooth dynamical systems (Buzzi et al. \cite{Bca-bsm}), perturbed Gardner equation (Wen \cite{wen-MMAS}), delayed Schr\"{o}dinger equation (Xu et al. \cite{XCH-JAAC}), perturbed mKdV and mK$(3,1)$ equation (Zhang et al. \cite{Zhanglijun-ND,Zhanglijun-IJBC}),  generalized Keller-Segel system (Du et al. \cite{DZ-JDE2}), and so on.

To apply the geometric singular perturbation approach to track invariant manifolds of corresponding ordinary differential equations (ODEs), usually, it is connected with the zeroes of the Melnikov functions associated to the perturbed ODEs. However, it is not easy to determine the zeroes of the  Melnikov functions because of the complexity of the  expression.
To analyze the Melnikov functions, usually, some auxiliary tools are needed. For examples, one of the effective method is to detect the monotonicity of the ratio of Abelian integral (for short, MRAI) associated to the Melnikov function.  For instance, Derks and Gils \cite{DG-JJIAM}, Ogawa \cite{OT-HMJ} computed the MRAI of perturbed KdV equation. Chen et al. \cite{CAY-JDE,CAY-AML,CAY-QTDS} detected the MRAI of perturbed BBM equation, perturbed defocusing mKdV equation and perturbed generalized KdV equation by the similar method. Du et al. \cite{DZJ-SSM} considered MRAI of a generalized Nizhnik-Novikov-Veselov equation with diffusion term. Then they all proved the existence of traveling wave solutions. Different from them, Sun et al. \cite{SunNA} employed the Chebyshev criterion to detect the MRAI of a shallow water fluid to analyze the coexistence of the solitary and periodic wave solutions.  Sun and Yu \cite{SY2019} also illustrated the existence and uniqueness of periodic waves of a generalized BBM equation based on the same technique.

Different from applying the method of MRAI (\cite{CAY-JDE,CAY-AML,CAY-QTDS,DZJ-SSM,DG-JJIAM,OT-HMJ,SunNA,SY2019}) and the method of Fredholm orthogonality (Du and Qiao \cite{DQ2020}), in this paper, we directly compute the explicit expression of Melnikov function to determine the zeroes of the Melnikov function.
The main purpose of this paper is to investigate the existence of solitary and kink (anti-kink) wave solutions for $(1+1)$-dimensional DLWE under a small singular perturbation. Firstly, we introduce a new parameter ${\kappa}$ (horizontal ordinate of the saddle point associated to the homoclinic orbit) to obtain the solitary wave solutions and kink (anti-kink) wave solutions of the unperturbed $(1+1)$-dimensional DLWE. Then, the geometric singular perturbation approach is used to construct homoclinic or heteroclinic orbits by tracking invariant manifolds of corresponding ODEs. Finally, by Melnikov method, we analyze conditions such as the wave speed and parameter for the existence of solitary and kink (anti-kink) wave solutions. The contributions and novelty of this paper is summarized as follows:\\
\noindent (1) By GSP approach and the bifurcation analysis, we prove that the solitary wave solution exists at a suitable wave speed $c$ for the bifurcation parameter ${\kappa}\in (0,1-\frac{\sqrt3}{3})\cup (1+\frac{\sqrt3}{3},2)$, while the kink and anti-kink wave solutions exist at a unique wave speed $c^*=\sqrt{15}/3$ for $\kappa=0$ or $\kappa=2$.  {This indicates that the solitary and kink (anti-kink) wave solutions can't coexist at the same bifurcation parameter value in $(1+1)$-dimensional DLWE.}\\
\noindent (2) To prove the existence of solitary and kink (anti-kink) wave solutions, we use the analytical Melnikov method and bifurcation theory. Not only the analytical expression of Melnikov integral is directly obtained by hand (not by mathematical software) for a perturbed PDE, but also the analytical expression of the limit wave speed $c$ is directly obtained by hand (see \eqref{Hom-Melnikov}-\eqref{c-Melnikov1} in Section 5). Indeed, the calculations and expressions are verified by the mathematical softwares (e.g., Maple).\\
\noindent (3) To obtained the analytical expressions, we use the following two mathematical skills.\\
\indent {\bf (i)} We apply a factorization technique together with the handbook of integral \cite{GI-SV} to the calculate $I(\kappa)$ and $J(\kappa)$, which lead us to obtain analytical Melnikov integral in section 5.
We perform very detailed computations and skillful mathematical analysis to obtain the expressions of $I(\kappa)$ and $J(\kappa)$  (see \eqref{I}--\eqref{J22} in Section 5).\\
\indent {\bf (ii)} $I{({\kappa})}$ and $J{({\kappa})}$ are integrated over a closed curve. Usually, it is difficult to be calculated. We transfer them to an definite integral by translating the time variable $\zeta$ into the state variable $z$ on homoclinic orbit. Then the up and lower limits of the definite integral depend on the coordinates of the intersections between the homoclinic orbit and $z$ axis. For example,
\[
I({\kappa})=\displaystyle\oint_{\Gamma({\kappa})}y^2d\zeta=\displaystyle\oint_{\Gamma({\kappa})}ydz
=2 \displaystyle\int_{z_+}^{{\kappa}}f(z) dz,
\]
where $f(z)$ is the expression related to the homoclinic orbit (see \eqref{I} in Section 5).

\noindent (4) We introduce a new bifurcation parameter ${\kappa}$ (horizontal ordinate of the saddle point associated to the homoclinic orbit) instead of the integral constant $G$. This benefits to the factorization of the integrand because there is a factor $(z-{\kappa})^2$. For example, when we compute $I({\kappa})$, it will lead us to
\[
I({\kappa})=\cdots=\displaystyle\int_{z_+}^{{\kappa}}\sqrt{\frac{3}{4}(z-{\kappa})^2(z-z_+)(z-z_-)}dz=\cdots
\]
(see \eqref{I} in section 5).

\noindent (5) Numerical simulations are utilized to verify the theoretical results.

\subsection{Outline of the paper}

The outline of the paper is as follows. Some preliminaries including geometric singular perturbation theory are introduced in section 2.  { Reduction of the model by geometric singular perturbation theory and dynamical system method is presented in Section 3.}
%Bifurcations of phase portraits for the unperturbed $(1+1)$-dimensional DLWE based on dynamical system method are obtained in section 3.
By introducing a bifurcation parameter ${\kappa}$, the solitary and kink (anti-kink) wave solutions are determined in section 4. In section 5, combining the GSP approach and Melnikov method, the existence of solitary and kink (anti-kink) wave solutions for a perturbed $(1+1)$-dimensional DLWE is investigated.  In section 6, numerical simulations are carried out to show the effectiveness of previous theoretical results. %Finally a conclusion ends our work.
\section{Preliminaries}
In this section, we introduce some known results on the theory of geometric singular perturbation (see. e.g., \cite{JCK-berlin}).
Consider the system
\begin{equation}\label{System1}
\begin{cases}
x_1'=f(x_1,x_2,\epsilon),\\
x_2'=\epsilon g(x_1,x_2,\epsilon),
\end{cases}
\end{equation}
where $'=\frac{d}{dt}$, $x_1\in \mathbb{R}^n$, $x_2\in \mathbb{R}^l$ and $\epsilon$ is a positive real parameter, $U\subseteq \mathbb{R}^{n+l}$ is open subset, and $I$ is an open subset of $\mathbb{R}$, containing $0$. $f$ and $g$ are $C^\infty$ on a set $U\times I$. Moreover, the $x_1$ (resp., $x_2$) variables are called fast (resp., slow) variables. Letting $\tau=\epsilon t$ which gives the following equivalent system
\begin{equation}\label{System2}
\begin{cases}
\epsilon \dot{x_1}=f(x_1,x_2,\epsilon),\\
\dot{x_2}= g(x_1,x_2,\epsilon),
\end{cases}
\end{equation}
where $\cdot=\frac{d}{d\tau}$. We refer to $t$ (resp., $\tau$ ) as the fast time scale or fast time (resp., slow time scale or slow time). Each of the scalings is naturally associated with a limit as $\epsilon$ tend to zero. These limits are respectively given by
\begin{equation}\label{System3}
\begin{cases}
x_1'=f(x_1,x_2,0),\\
x_2'=0,
\end{cases}
\end{equation}
and
\begin{equation}\label{System4}
\begin{cases}
0=f(x_1,x_2,0),\\
x_2'=g(x_1,x_2,0).
\end{cases}
\end{equation}
System (\ref{System3}) is called the layer problem and system (\ref{System4}) is reduced system.
\begin{mydef}\label{definition1}
(see \cite{JCK-berlin,DLL2018,cheng-dcds}) A manifold $M_0$ on which $f(x_1,x_2,0)=0$ is called a critical manifold or slow manifold. A critical manifold $M_0$ is said to be normally hyperbolic if the linearization of system (\ref{System1}) at each point in $M_0$ has exactly $l$ eigenvalues on the imaginary axis $Re(\lambda)=0$.
\end{mydef}
\begin{mydef}\label{definition2}
(see \cite{JCK-berlin,DLL2018,cheng-dcds}) A set $M$ is locally invariant under the flow of system (\ref{System1}) if it has neighborhood $V$ so that no trajectory can leave $M$ without also leaving $V$. In other words, it is locally invariant if for all $x_1\in M$, $x_1\cdot [0,t]\subseteq V$ implies that $x_1\in M$, $x_1\cdot [0,t]\subseteq M$, similarly with $[0,t]$ replaced by $[t,0]$ when $t<0$, where $x_1\cdot [0,t]$ denotes the application of a flow after time $t$ to the initial condition $x_1$.
\end{mydef}
\begin{mylem}\label{Lemma1.1}
(see \cite{FN-JDE,DLL2018,cheng-dcds}) Let $M_0$ be a compact, normally hyperbolic critical manifold given as a graph $\{(x_1,x_2):x_1=h^0(x_2)\}$. Then for sufficiently small positive $\epsilon$ and any $0<r<+\infty$,

$\bullet$ there exists a manifold $M_\epsilon$, which is locally invariant under the flow of system (\ref{System1}) and $C^r$ in $x_1$, $x_2$, $\epsilon$. Moreover, $M_\epsilon$ is given as graph:
$$M_\epsilon=\{(x_1,x_2):x_1=h^\epsilon(x_2)\}$$
for some $C^r$ function $h^\epsilon(x_2)$;

$\bullet$ $M_\epsilon$ possesses locally invariant stable and unstable manifold $W^s(M_\epsilon)$ and $W^u(M_\epsilon)$ lying within $O(\epsilon)$ and being $C^r$ diffeomorphic to the stable and unstable manifold $W^s(M_0)$ and $W^u(M_0)$ of the critical manifold $M_0$;

$\bullet$ $W^s(M_\epsilon)$ is partitioned by moving invariant submanifolds $\digamma^s(p_\epsilon)$, which are $O(\epsilon)$ close and diffeomorphic to $\digamma^s(p_0)$, with base point $p_{\epsilon}$ belonging to $M_{\epsilon}$. Moreover, they are $C^r$ with respect to $p$ and $\epsilon$. Moving invariance means the submanifold $\digamma^s(p_\epsilon)$ is mapped under the time $t$ flow to another submanifold $\digamma^s(p_\epsilon\cdot t)$ whose base point is the time $t$ evolution image of the taken base point $p_\epsilon$;

$\bullet$ the dynamics on $M_\epsilon$ is a regular perturbation of that generated by system (\ref{System4}).
\end{mylem}

\section{Reduction of the model by geometric singular perturbation theory and dynamical system method  }
In this section,  to consider the exact solutions of the unperturbed $(1+1)$-dimensional DLWE (\ref{wuzhang1}), we firstly reduce the model by GSP theory and dynamical system method. From the view of physical meanings, we only consider the case of $c>0$ in this paper.

By introducing the following transformations:
\begin{equation}\label{transformations1}
u(x,t)=\phi(\xi),\ \ \ v(x,t)=\psi(\xi),\ \ \ \xi=x-ct,
\end{equation}
we obtain
\begin{equation}\label{transformations2}
\begin{array}{l}
\frac{\partial u(x,t)}{\partial t}=-c\phi_{\xi},\ \  \frac{\partial v(x,t)}{\partial t}=-c\psi_{\xi},\\
\frac{\partial^{2 }u(x,t)}{\partial x^{2 }}=\phi_{\xi\xi},\ \ \frac{\partial^{3 }u(x,t)}{\partial x^{3 }}=\phi_{\xi\xi\xi},\ \  \frac{\partial^{4 }u(x,t)}{\partial x^{4 }}=\phi_{\xi\xi\xi\xi},
\end{array}
\end{equation}
where $ \phi_{\xi}$ and $ \psi_{\xi}$ are the first order derivative with respect to $\xi$. And $ \phi_{\xi\xi}$, $ \phi_{\xi\xi\xi}$ and $ \phi_{\xi\xi\xi\xi}$ means the second, third and fourth order derivative with respect to $\xi$, respectively.

Substituting (\ref{transformations1}) and (\ref{transformations2}) into (\ref{wuzhang1}), then system (\ref{wuzhang1}) is given by
\begin{equation}\label{ODEs}
\begin{cases}
c \phi_{\xi}=\phi \phi_{\xi}+\psi_{\xi},\\
c \psi_{\xi}=\psi \phi_{\xi}+\phi \psi_{\xi}+\frac{1}{3}\phi_{\xi\xi\xi}-\varepsilon (\phi_{\xi\xi}+\phi_{\xi\xi\xi\xi}).
\end{cases}
\end{equation}
Integrating both sides of the first equation of system (\ref{ODEs}) once with respect to $\xi$  and letting the integration constant be zero yields
\begin{equation}\label{Integration}
\psi(\xi)=c \phi(\xi)-\frac{1}{2}\phi^2(\xi).
\end{equation}
We substitute (\ref{Integration}) into the second equation of (\ref{ODEs}), then the coupled system (\ref{ODEs}) becomes the following ODE:
\begin{equation}\label{ODE1}
\frac{1}{3}\phi_{\xi\xi\xi}-\frac{3}{2}\phi^2\phi_{\xi}+3c \phi \phi_\xi-c^2\phi_\xi-\varepsilon (\phi_{\xi\xi}+ \phi_{\xi\xi\xi\xi})=0.
\end{equation}
Integrating both sides on (\ref{ODE1}) once with respect to $\xi$ and rescaling $\epsilon=3\varepsilon$, it follows that
\begin{equation}\label{ODE2}
\phi_{\xi\xi}-\frac{3}{2}\phi^3+\frac{9}{2}c\phi^2-3c^2\phi-\epsilon(\phi_{\xi}+ \phi_{\xi\xi\xi})=g,
\end{equation}
where $g$ is an integration constant ($ g\in \mathbb{R} $).

Introducing new variables $\zeta=c\xi$, $z=\frac{\phi}{c}$ and $G=\frac{g}{c^3}$, Eq. (\ref{ODE2}) is equivalent to
\begin{equation}\label{ODE3}
-3z+\frac{9}{2}z^2-\frac{3}{2}z^3+\frac{d^2z}{d\zeta^2}-\epsilon(\frac{1}{c}\frac{dz}{d\zeta}+c\frac{d^3z}{d\zeta^3})=G.
\end{equation}
Obviously, Eq. (\ref{ODE3}) reduces to a three-dimensional system:
\begin{equation}\label{ODE4}
\begin{cases}
\frac{dz}{d\zeta}=y,\\
\frac{dy}{d\zeta}=w,\\
\epsilon c\frac{dw}{d\zeta}=-3z+\frac{9}{2}z^2-\frac{3}{2}z^3+w-\frac{\epsilon}{c}y-G.
\end{cases}
\end{equation}
where $\epsilon$ is a sufficiently small parameter such that $0<\epsilon\ll1$. Therefore, the traveling wave solutions of Eq. (\ref{wuzhang1}) can be obtained by studying the corresponding orbits of system (\ref{ODE4}).

Obviously, system (\ref{ODE4}) is a singularly perturbed system described as ``slow system".  Rescaling $\zeta=\epsilon \eta $, we have the following equivalent ``fast system":
\begin{equation}\label{ODE4b}
\begin{cases}
\frac{dz}{d\eta}=\epsilon y,\\
\frac{dy}{d\eta}=\epsilon w,\\
 c\frac{dw}{d\eta}=-3z+\frac{9}{2}z^2-\frac{3}{2}z^3+w-\frac{\epsilon}{c}y-G.
\end{cases}
\end{equation}
Let $\epsilon=0$ in system (\ref{ODE4}) and  (\ref{ODE4b}). Then, the corresponding reduced system is given by
\begin{equation}\label{ODE4c}
\begin{cases}
\frac{dz}{d\zeta}=y,\\
\frac{dy}{d\zeta}=w,\\
0=-3z+\frac{9}{2}z^2-\frac{3}{2}z^3+w-G,
\end{cases}
\end{equation}
and the layer system is
\begin{equation}\label{ODE4d}
\begin{cases}
\frac{dz}{d\eta}=0,\\
\frac{dy}{d\eta}=0,\\
 c\frac{dw}{d\eta}=-3z+\frac{9}{2}z^2-\frac{3}{2}z^3+w-G.
\end{cases}
\end{equation}
which admits a two-dimensional critical manifold as follows
\begin{equation}\label{hyperbolic}
M_0=\{(z,y,w) \in \mathbb{R}^3|w=3z-\frac{9}{2}z^2+\frac{3}{2}z^3+G \}.
\end{equation}
Suppose $A$ to be the linearized matrix of system  (\ref{ODE4d}), then it is given by
\begin{equation}
A=\left(
\begin{array}{ccc}
0 & 0 & 0 \\
0 & 0 & 0 \\
\frac{-3+9z-\frac{9}{2}z^2}{c} & 0 & \frac{1}{c} \\
\end{array}
\right),
\end{equation}
that the eigenvalues of A are $0$, $0$ and $\frac{1}{c}$. Thus, $M_0$ is a normally hyperbolic invariant manifold (see \textbf{Definition} \ref{definition1} and \textbf{Definition} \ref{definition2}). There exists a two-dimensional slow submanifold $M_\epsilon$ which is $C^1$ $O(\epsilon)$ close to $M_0$ (see  \textbf{Lemma} \ref{Lemma1.1}). The invariant submanifold $M_\epsilon$ is represented by
\begin{equation}\label{hyperbolic}
M_\epsilon=\{(z,y,w) \in \mathbb{R}^3|w=3z-\frac{9}{2}z^2+\frac{3}{2}z^3+G+\epsilon[cy(\frac{9}{2}z^2-9z+3+\frac{1}{c^2})]+O(\epsilon^2) \}.
\end{equation}
The dynamical behaviour of slow system (\ref{ODE4}) or fast system (\ref{ODE4b}) is governed by
\begin{equation}\label{ODE5}
\begin{cases}
\frac{dz}{d\zeta}=y,\\
\frac{dy}{d\zeta}=3z-\frac{9}{2}z^2+\frac{3}{2}z^3+G+\epsilon[cy(\frac{9}{2}z^2-9z+3+\frac{1}{c^2})]+O(\epsilon^2).
\end{cases}
\end{equation}
Obviously, the unperturbed system $(\ref{ODE5})\mid_{\epsilon=0}$ admits homoclinic, heteroclinic and periodic orbits with the parameter $G$ taking different values. By the bifurcation theory of planar dynamical systems (see \cite{Li-Book}), we have the following proposition of unperturbed system $(\ref{ODE5})\mid_{\epsilon=0}$:
\begin{myprop}\label{proposition}
(i) If $\mid G\mid>\frac{\sqrt3}{3}$, there exists a single saddle point (see Fig.1 (a) or (g)).\\
(ii) If $\mid G\mid=\frac{\sqrt3}{3}$, there exist a saddle point and a cusp point (see Fig.1 (b) or (f)).\\
(iii) If $G=0$, there exist two saddle points $(0,0)$ and $(2,0)$, and a center point $(1,0)$. It has two heteroclinic orbits surrounding the center point $(1,0)$ to two saddle points $(0,0)$ and $(2,0)$ (see Fig.1 (d)).\\
(iv) If $0<\mid G\mid<\frac{\sqrt3}{3}$, there exist two saddle points and a center point. It has a homoclinic orbit surrounding the center point to one saddle point (see Fig.1 (c) or (e)).\end{myprop}
\begin{center}
\begin{tabular}{cccc}
\epsfxsize=3.2cm \epsfysize=3.2cm \epsffile{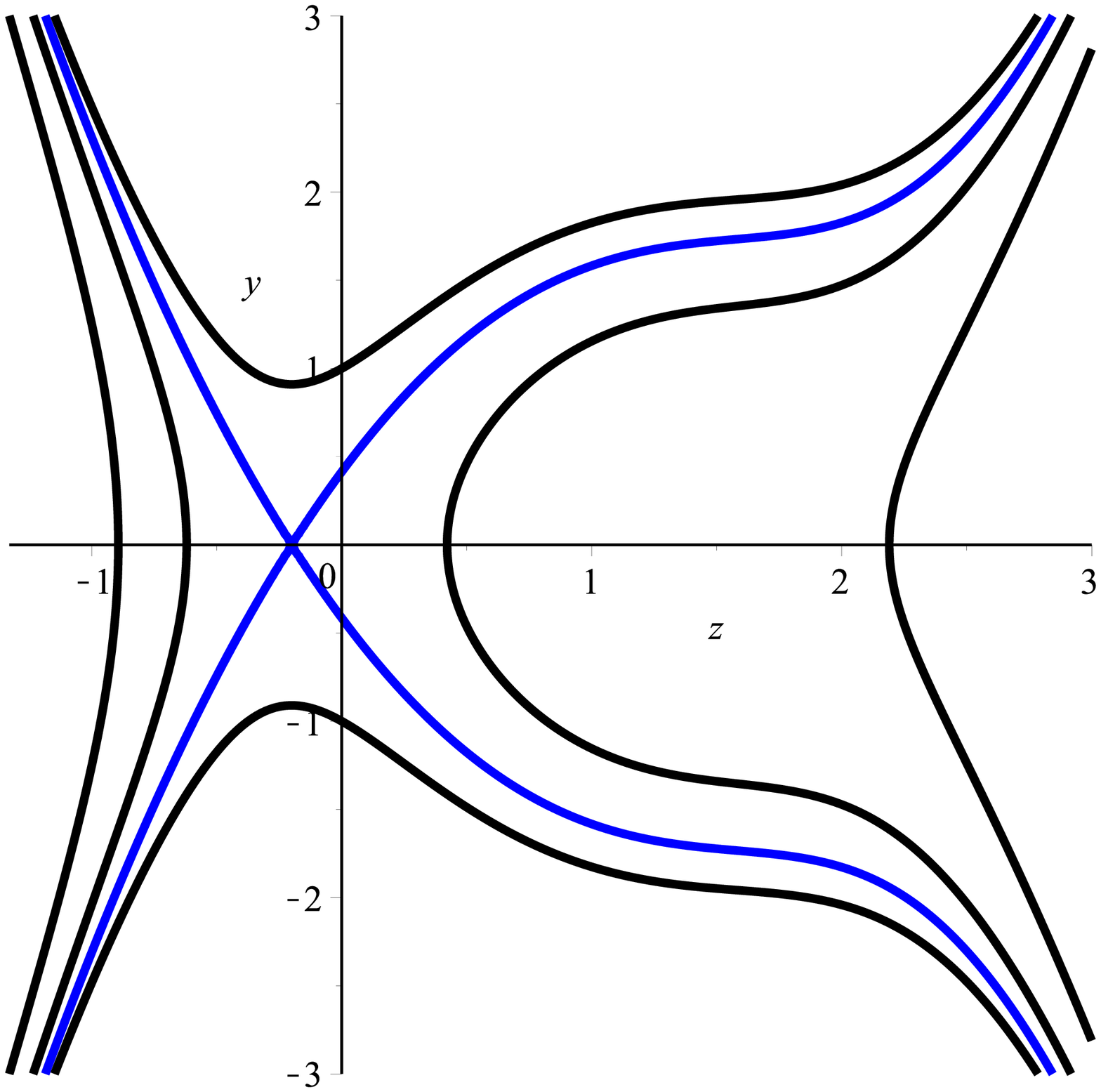}&
\epsfxsize=3.2cm \epsfysize=3.2cm \epsffile{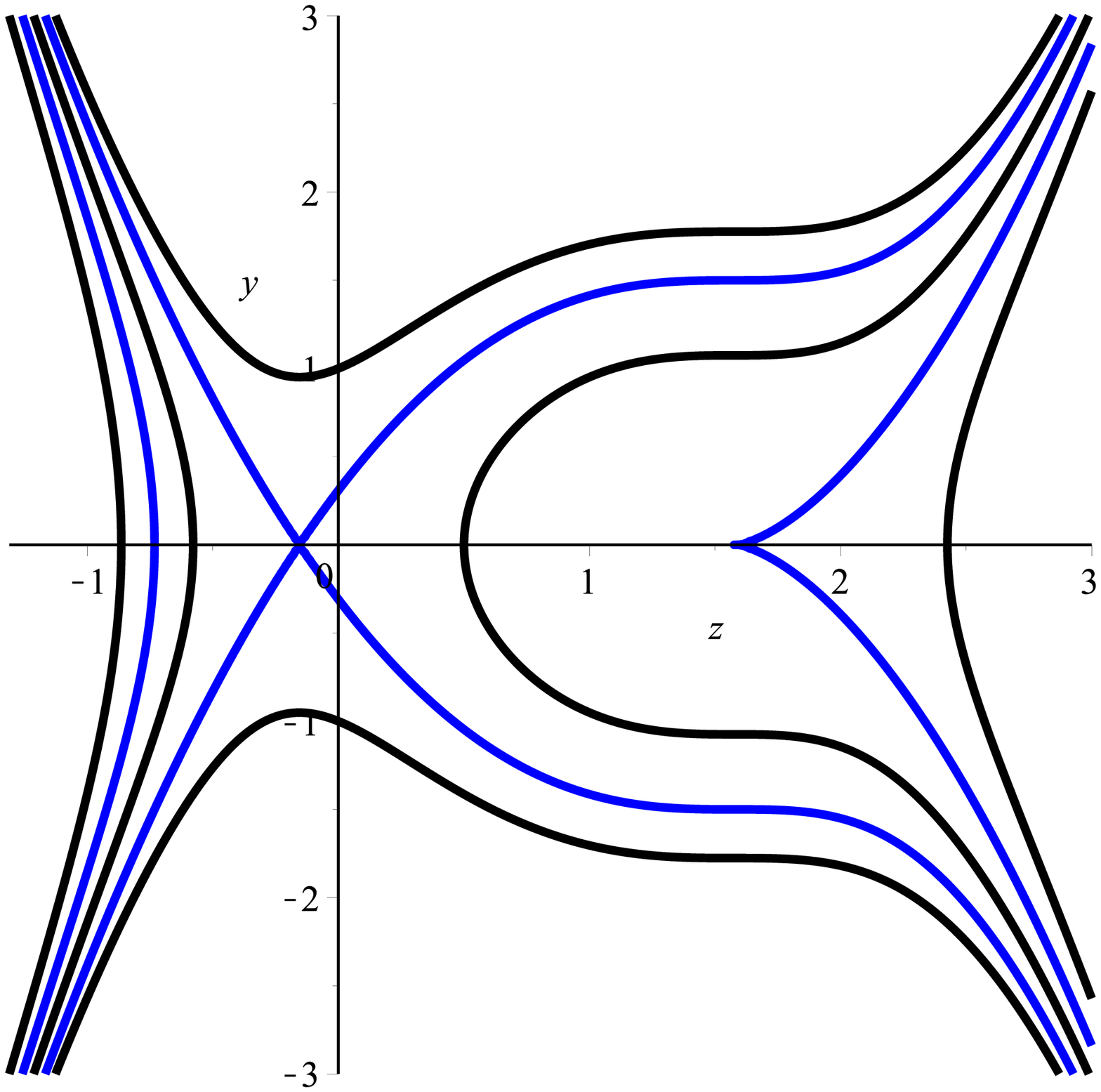}&
\epsfxsize=3.2cm \epsfysize=3.2cm \epsffile{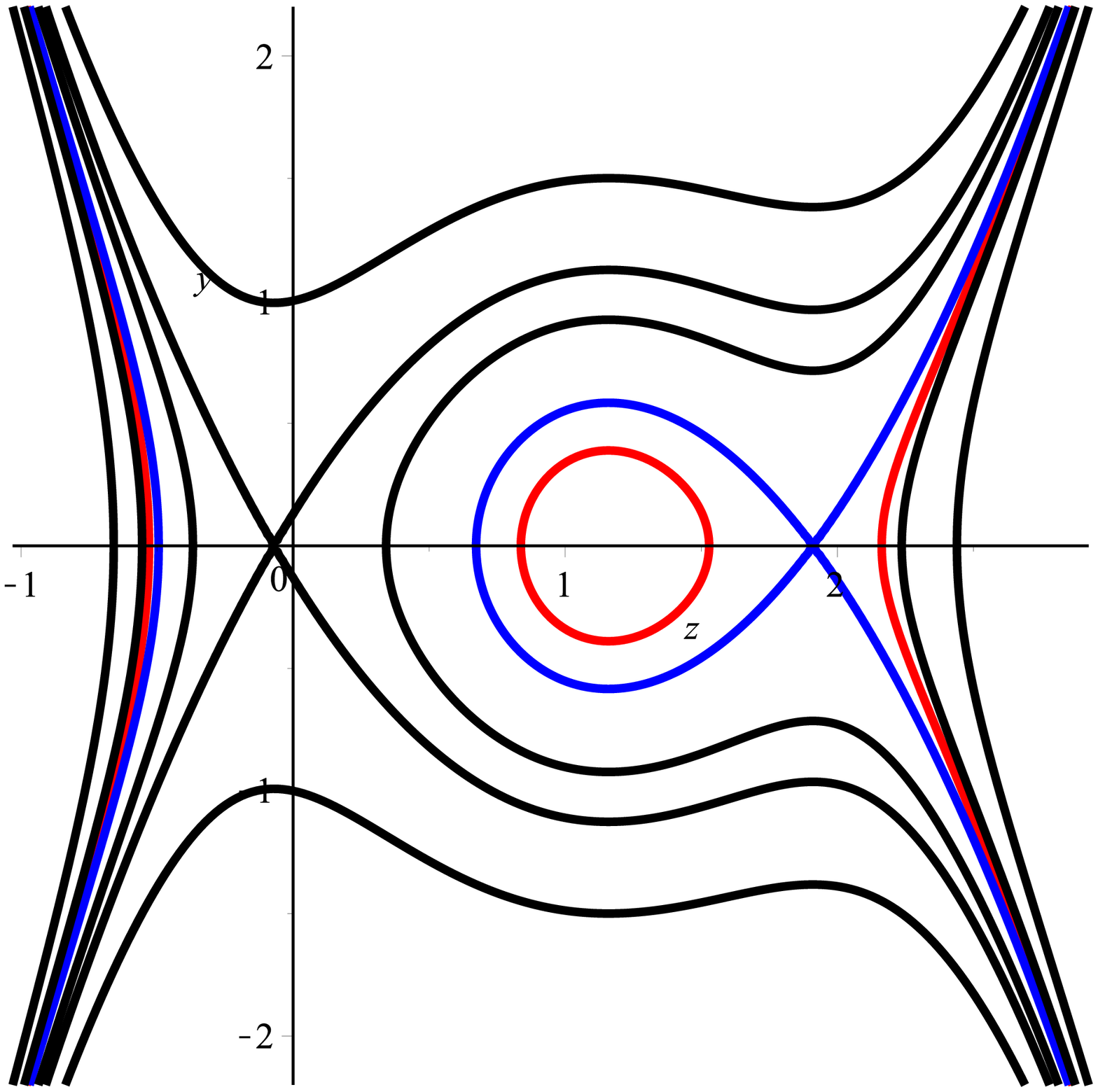}&
\epsfxsize=3.2cm \epsfysize=3.2cm \epsffile{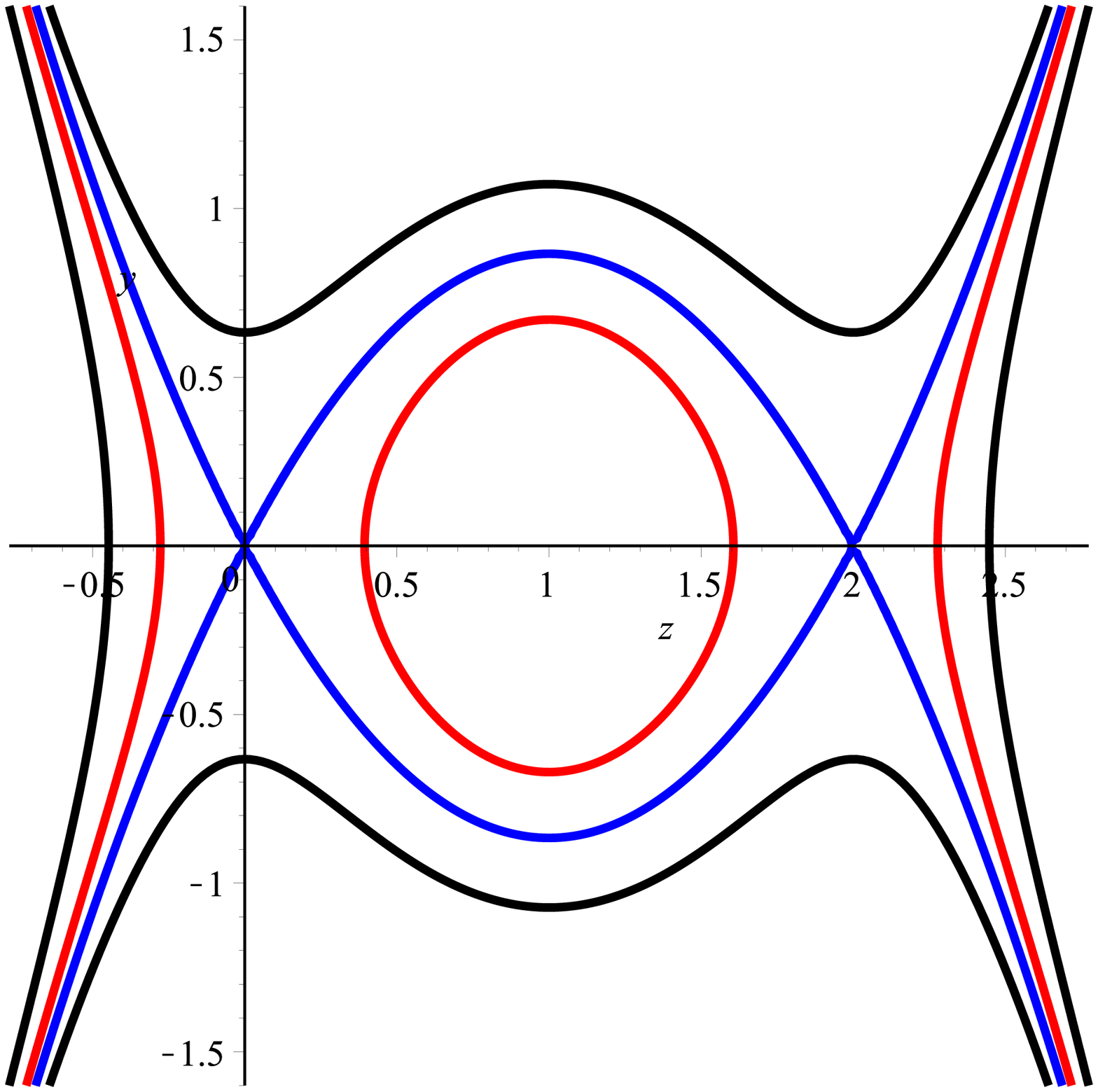}\\
\footnotesize{(a)\ \ $G<-\frac{\sqrt3}{3}$} & \footnotesize{ (b)\ \ $G=-\frac{\sqrt3}{3}$ } & \footnotesize{ (c)\ \ $-\frac{\sqrt3}{3}<G<0$} & \footnotesize{ (d)\ \ $G=0$}
\end{tabular}
\end{center}
\begin{center}
\begin{tabular}{ccc}
\epsfxsize=3.2cm \epsfysize=3.2cm \epsffile{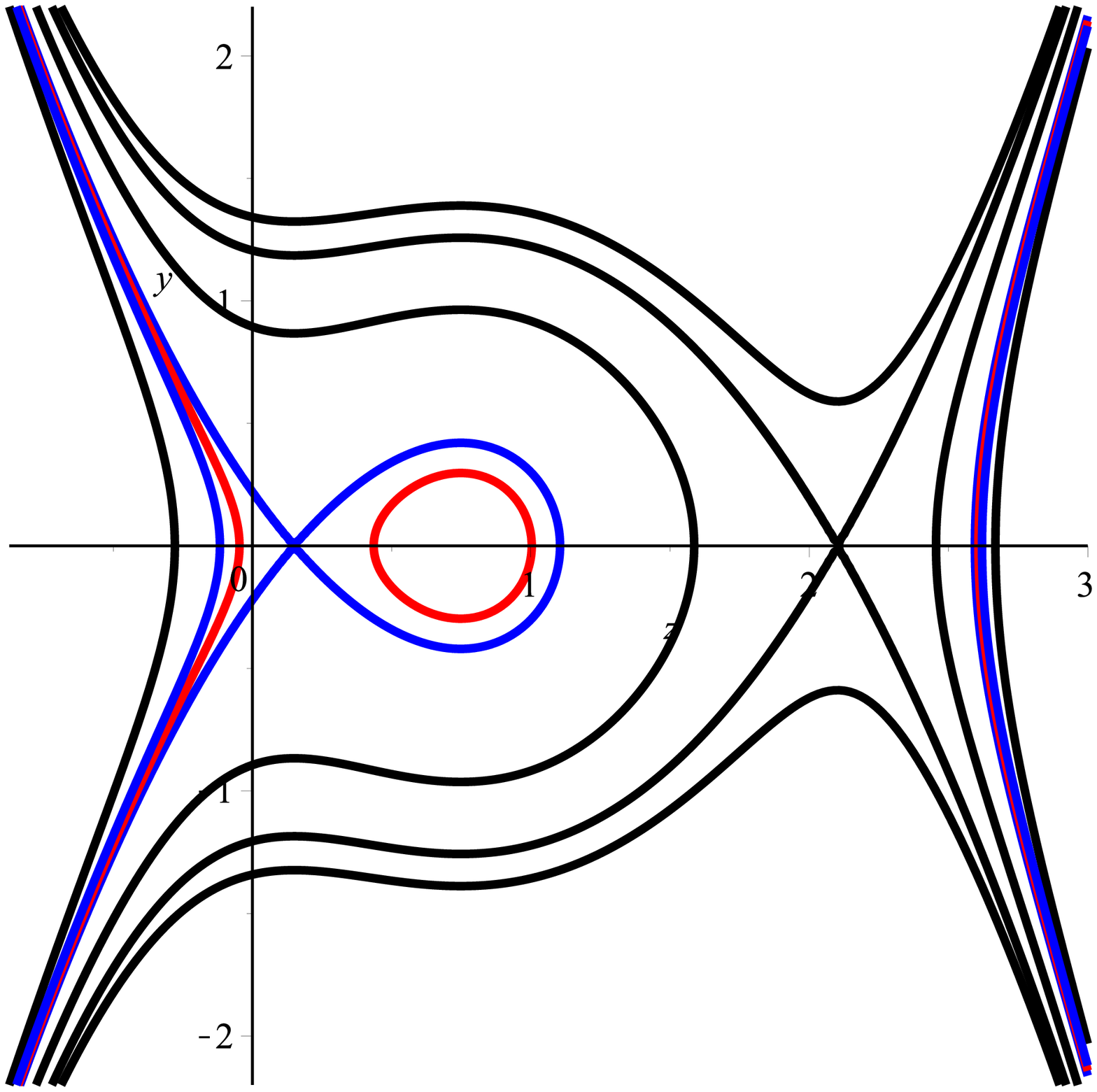}&
\epsfxsize=3.2cm \epsfysize=3.2cm \epsffile{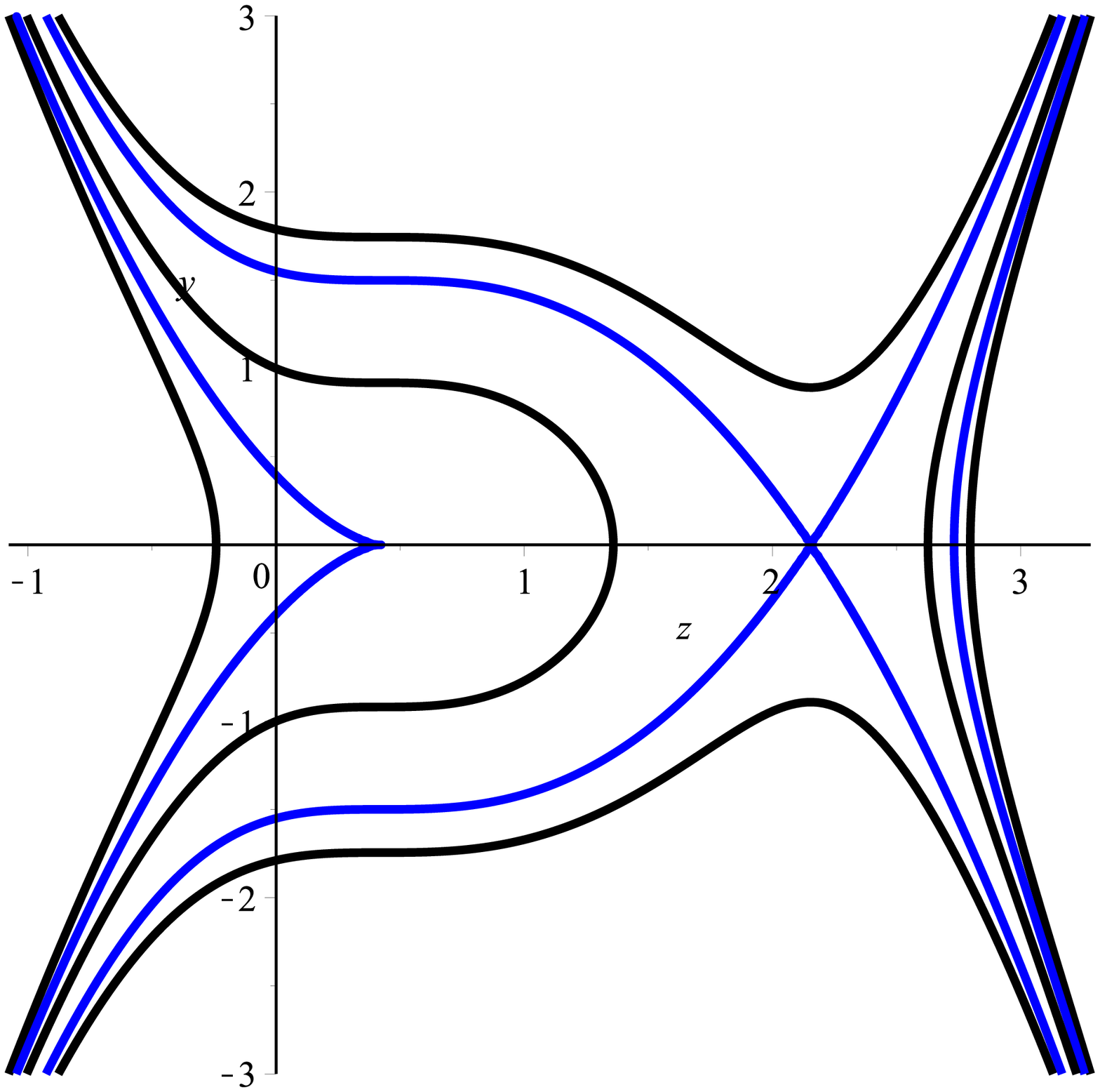}&
\epsfxsize=3.2cm \epsfysize=3.2cm \epsffile{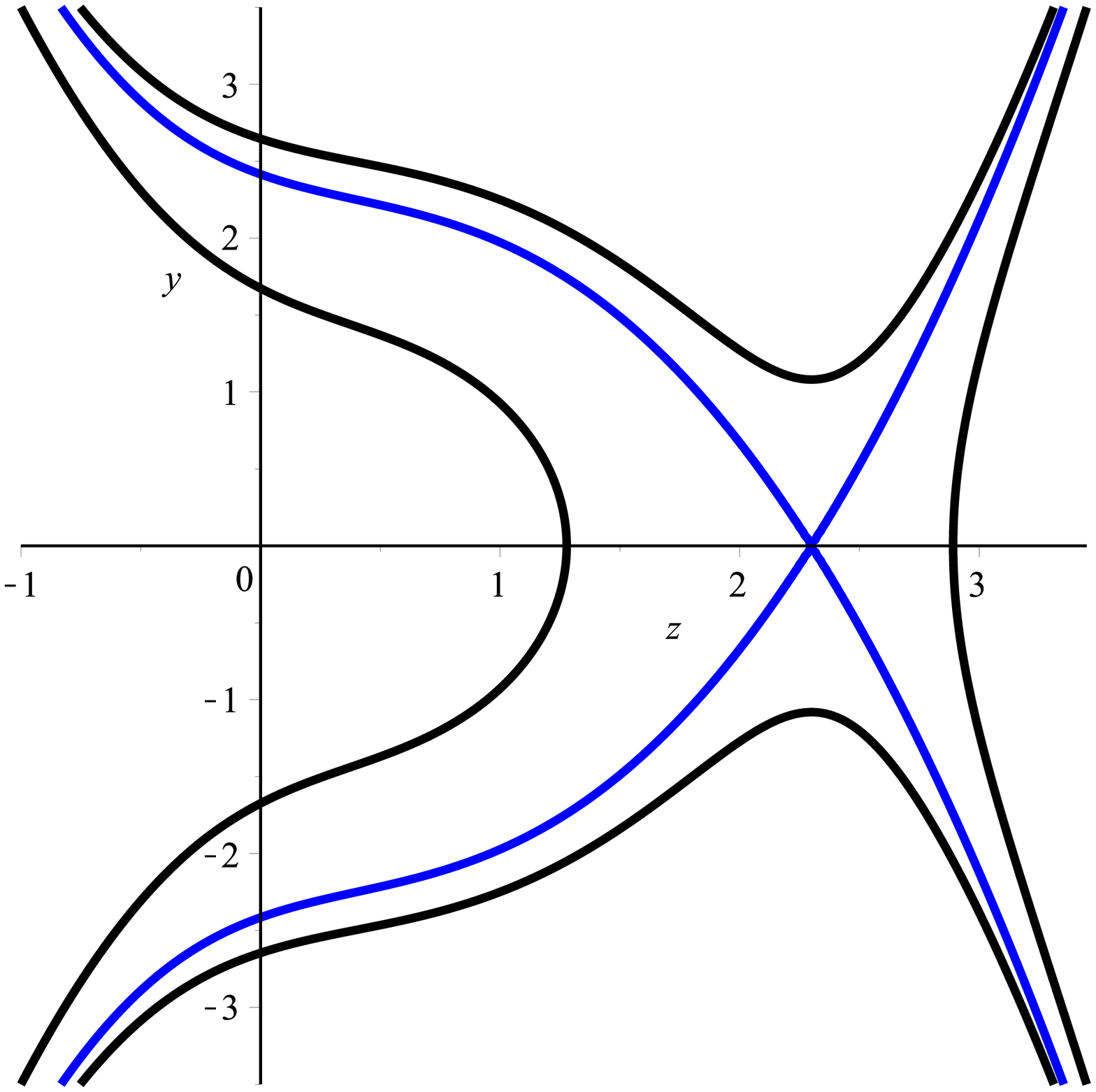}\\
\footnotesize{(e)\ \ $0<G<\frac{\sqrt3}{3}$} & \footnotesize{ (f)\ \ $G=\frac{\sqrt3}{3}$} & \footnotesize{ (g) \ \ $G> \frac{\sqrt3}{3}$}
\end{tabular}
\end{center}
\begin{center}
\footnotesize {{Fig. 1} \ \ The bifurcation and phase portraits of unperturbed system (\ref{ODE5})$\mid_{\epsilon=0}$.}
\end{center}
Assume that $z(\zeta)$ is a solution of Eq. (\ref{ODE3}) satisfying $\lim\limits_{\zeta\longrightarrow \infty}z(\zeta)={\kappa}$, thus $u(x,t)=cz(\zeta)=cz[c(x-ct)]$ is a solitary wave solution or a kink wave solution of Eq. (\ref{ODE3}).

For $-\frac{\sqrt3}{3}<G\leq0$ and $0\leq G<\frac{\sqrt3}{3}$, the unperturbed system (\ref{ODE5})$\mid_{\epsilon=0}$ admits a saddle point defined by $S_0({\kappa},0)$. We know that $G=-(3{\kappa}-\frac{9}{2}{{\kappa}}^2+\frac{3}{2}{{\kappa}}^3)$, then we obtain $1+\frac{\sqrt3}{3}<{\kappa}\leq2$  for $G\in(-\frac{\sqrt3}{3},0]$ and $0\leq {\kappa}<1-\frac{\sqrt3}{3}$ for $G\in[0,\frac{\sqrt3}{3})$, respectively. We can regard ${\kappa}$ as a bifurcation parameter to find the homoclinic orbits of system (\ref{ODE5}). Therefore, we rewrite system (\ref{ODE5}), that is:
\begin{equation}\label{ODE6}
\begin{cases}
\frac{dz}{d\zeta}=y,\\
\frac{dy}{d\zeta}=3z-\frac{9}{2}z^2+\frac{3}{2}z^3-(3{\kappa}-\frac{9}{2}{{\kappa}}^2+\frac{3}{2}{{\kappa}}^3)+\epsilon[cy(\frac{9}{2}z^2-9z+3+\frac{1}{c^2})]+O(\epsilon^2).
\end{cases}
\end{equation}

\section{Traveling wave solutions of Eq. $(\ref{wuzhang1})\mid_{\varepsilon=0}$}
In this section, we study the solitary and kink (anti-kink) wave solutions of Eq. $(\ref{wuzhang1})\mid_{\varepsilon=0}$. By dynamical method, we consider the traveling wave solution of system $(\ref{ODE6})\mid_{\epsilon=0}$ with ${\kappa}\in [0,1-\frac{\sqrt3}{3})\cup (1+\frac{\sqrt3}{3},2]$. The unperturbed system $(\ref{ODE6})$ is of the form:
\begin{equation}\label{ODE7}
\begin{cases}
\frac{dz}{d\zeta}=y,\\
\frac{dy}{d\zeta}=3z-\frac{9}{2}z^2+\frac{3}{2}z^3-(3{\kappa}-\frac{9}{2}{{\kappa}}^2+\frac{3}{2}{{\kappa}}^3).
\end{cases}
\end{equation}
It is easy to obtain the first integral of system (\ref{ODE7}), it is given by
\begin{equation}\label{First integral}
H(z,y)=\frac{1}{2}y^2-\frac{3}{8}z^4+\frac{3}{2}z^3-\frac{3}{2}z^2+(\frac{3}{2}{{\kappa}}^3-\frac{9}{2}{{\kappa}}^2+3{\kappa})z.
\end{equation}
The homoclinic orbit $\Gamma({\kappa})$ to the saddle $({\kappa},0)$ is defined by
\begin{equation}\label{Hom}
\begin{array}{rcl}
H(z,y)&=&\frac{1}{2}y^2-\frac{3}{8}z^4+\frac{3}{2}z^3-\frac{3}{2}z^2+(\frac{3}{2}{{\kappa}}^3-\frac{9}{2}{{\kappa}}^2+3{\kappa})z\\
&=&\frac{9}{8}{{\kappa}}^4-3{{\kappa}}^3+\frac{3}{2}{{\kappa}}^2.
\end{array}
\end{equation}
And the heteroclinic orbits $\Upsilon({\kappa})_{\pm}$ (where $\pm$ represents upper and lower branch curves) to the two saddle points $(0,0)$ and $(2,0)$, namely, ${\kappa}=0$ or ${\kappa}=2$, the heteroclinic curve is determined by
\begin{equation}\label{Het}
\begin{array}{rcl}
H(z,y)=\frac{1}{2}y^2-\frac{3}{8}z^4+\frac{3}{2}z^3-\frac{3}{2}z^2=0.
\end{array}
\end{equation}
\subsection{Solitary wave solutions of Eq. $(\ref{wuzhang1})\mid_{\varepsilon=0}$}
By (\ref{Hom}) and \textbf{Proposition} \ref{proposition} (\emph{iv}), we can obtain the expression of $y$ as follows:
\begin{equation}\label{expression1}
\begin{array}{rcl}
y&=&\pm \sqrt{\frac{3}{4}z^4-3z^3+3z^2-(3{{\kappa}}^3-9{{\kappa}}^2+6{\kappa})z+\frac{9}{4}{{\kappa}}^4-6{{\kappa}}^3+3{{\kappa}}^2}\\
 &=&\pm \sqrt{\frac{3}{4}(z-{\kappa})^2(z-z_+)(z-z_-)},
\end{array}
\end{equation}
where $z_{\pm}=-{\kappa}+2\pm \sqrt{-2{{\kappa}}^2+4{\kappa}}$. It can be seen that $z_-<z_+<{\kappa}$ for $1+\frac{\sqrt3}{3}<{\kappa}<2$ (see Fig. 1(c)) and ${\kappa}<z_-<z_+$ for $0<{\kappa}<1-\frac{\sqrt3}{3}$ (see Fig. 1(e)).

Since $\frac{dz}{d\zeta}=y$, it implies
\begin{equation}\label{expression2}
\frac{dz}{d\zeta}=\pm \sqrt{\frac{3}{4}(z-{\kappa})^2(z-z_+)(z-z_-)},
\end{equation}
which yields
\begin{equation}\label{expression3}
\zeta=\pm \int_{z_+}^{z}\frac{1}{\sqrt{\frac{3}{4}(s-{\kappa})^2(s-z_+)(s-z_-)}}ds,\ \ (1+\frac{\sqrt3}{3}<{\kappa}<2),
\end{equation}
and
\begin{equation}\label{expression4}
\zeta=\pm \int_{z_-}^{z}\frac{1}{\sqrt{\frac{3}{4}(s-{\kappa})^2(z_+-s)(z_--s)}}ds, \ \ (0<{\kappa}<1-\frac{\sqrt3}{3}).
\end{equation}
Thus, for $1+\frac{\sqrt3}{3}<{\kappa}<2$, the parametric representation of the dark solitary wave of system $(\ref{ODE6})|_{\epsilon=0}$ can be obtained by (see Fig. 2(a))
\begin{equation}\label{homoclinic1}
z(\zeta)=-\frac{2(6{{\kappa}}^2-12{\kappa}+4)}{2\textmd{cosh}(\frac{1}{2}\sqrt{18{{\kappa}}^2-36 {\kappa}+12}\zeta)\sqrt{-2{{\kappa}}^2+4{\kappa}}+4{\kappa}-4}.
\end{equation}
Then, it corresponds to a dark solitary wave solution (see Fig. 2(b)) of Eq. $(\ref{wuzhang1})\mid_{\varepsilon=0}$ is given by:
\begin{equation}\label{homoclinic2}
u(x,t)=-\frac{2(6{{\kappa}}^2-12{\kappa}+4)c}{2\textmd{cosh}[\frac{1}{2}\sqrt{18{{\kappa}}^2-36 {\kappa}+12}c(x-ct)]\sqrt{-2{{\kappa}}^2+4{\kappa}}+4{\kappa}-4}.
\end{equation}
\begin{center}
\begin{tabular}{cc}
\epsfxsize=6.5cm \epsfysize=6.5cm \epsffile{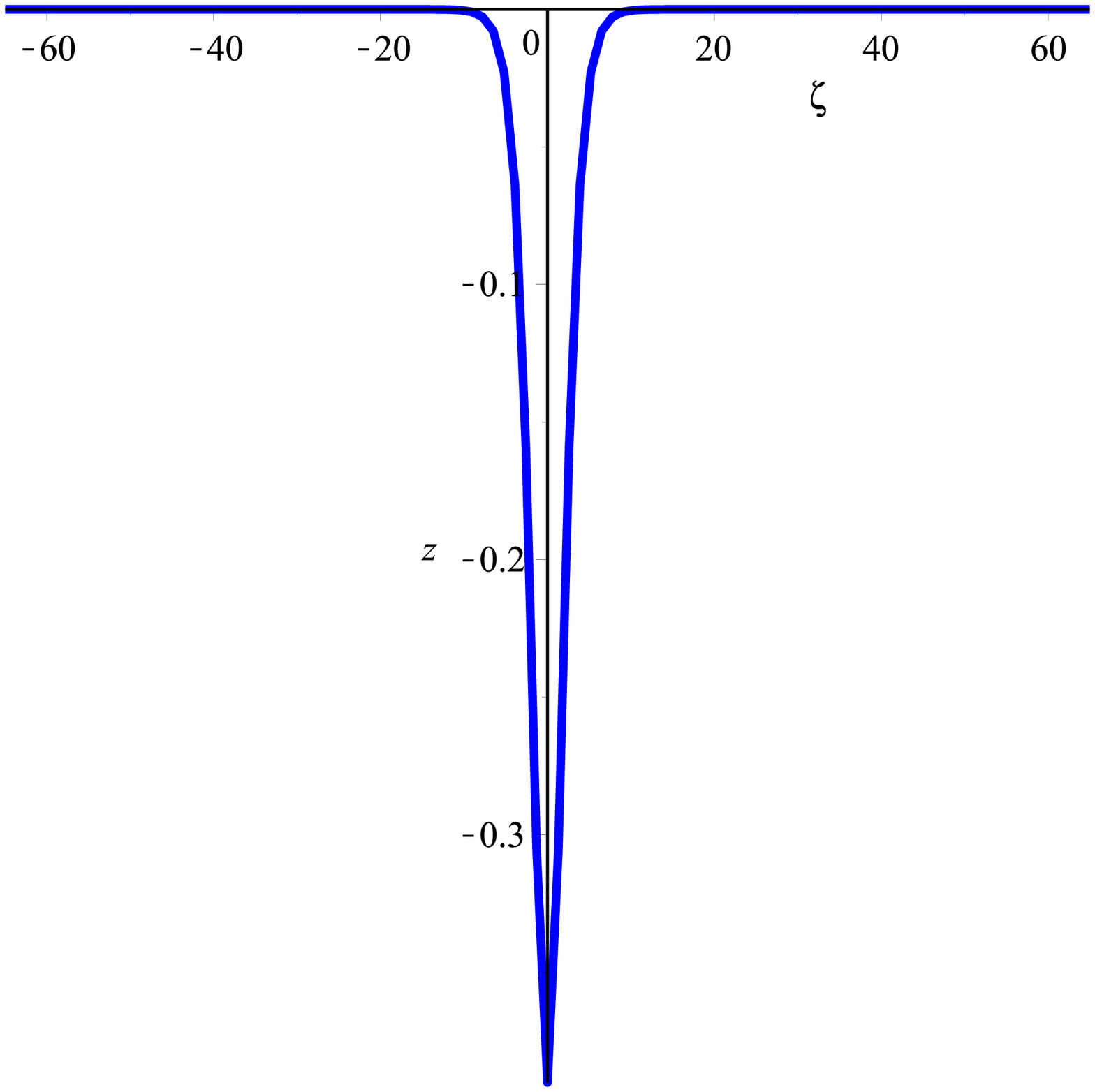}&
\epsfxsize=6.5cm \epsfysize=6.5cm \epsffile{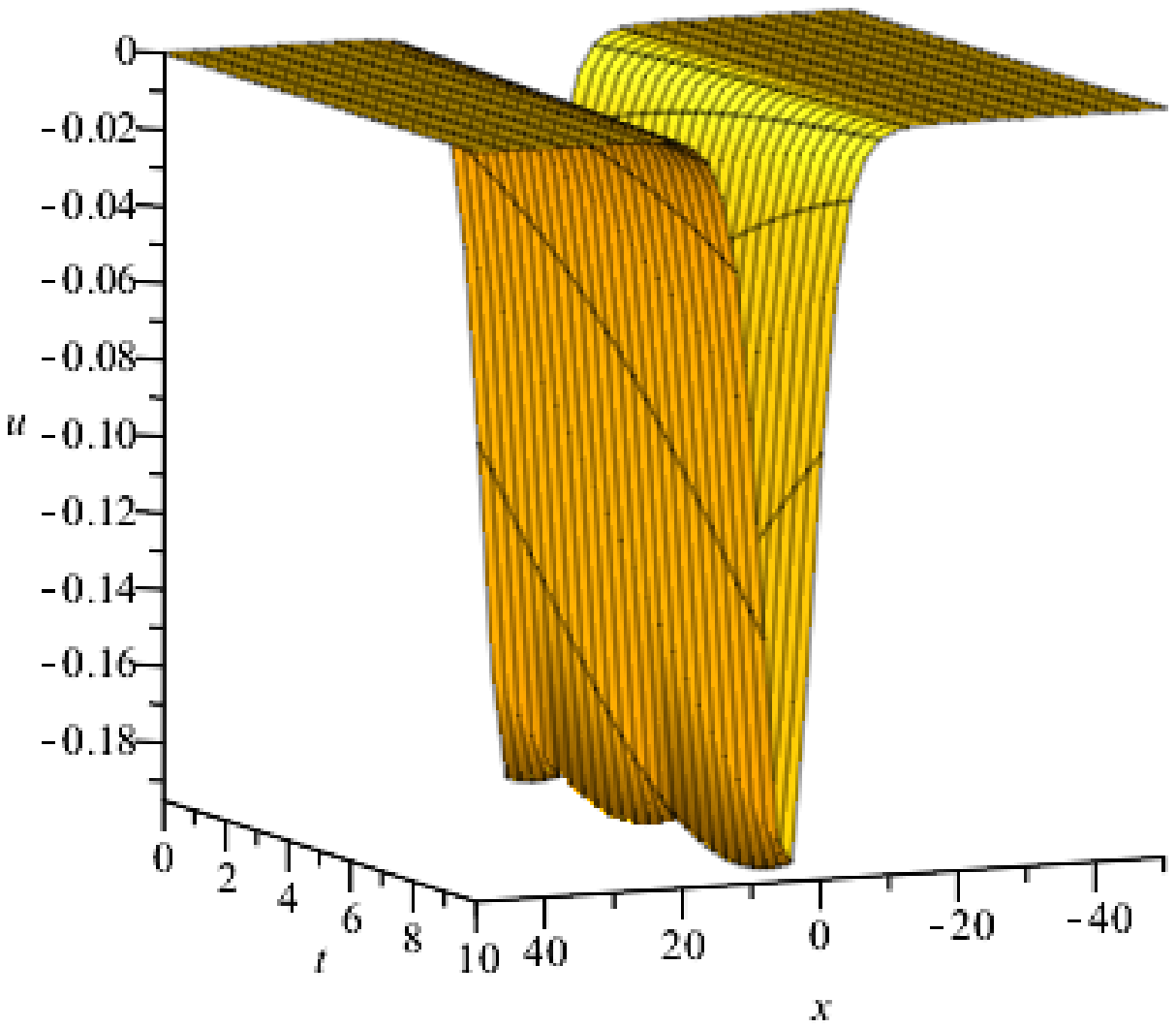}\\
\footnotesize{(a) Dark solitary wave } & \footnotesize{ (b) Exact solution  }
\end{tabular}
\end{center}
\begin{center}
\footnotesize {{Fig. 2} \ \ Dark solitary wave of system $(\ref{ODE6})|_{\epsilon=0}$ and exact solution of Eq. $(\ref{wuzhang1})\mid_{\varepsilon=0}$ for $c=0.5$ and ${\kappa}=1.7$ .}
\end{center}
For $0<{\kappa}<1-\frac{\sqrt3}{3}$, the expression of the bright solitary wave of system $(\ref{ODE6})|_{\epsilon=0}$ can be obtained by  (see Fig 3.(a))
\begin{equation}\label{homoclinic3}
z(\zeta)=-\frac{2(6{{\kappa}}^2-12{\kappa}+4)}{4{\kappa}-4-2\textmd{cosh}(\frac{1}{2}\sqrt{18{{\kappa}}^2-36 {\kappa}+12}\zeta)\sqrt{-2{{\kappa}}^2+4{\kappa}}}.
\end{equation}
The bright solitary wave solution of Eq. $(\ref{wuzhang1})\mid_{\varepsilon=0}$ (see Fig 3.(b)) is of the form
\begin{equation}\label{homoclinic4}
u(x,t)=-\frac{2(6{{\kappa}}^2-12{\kappa}+4)c}{4{\kappa}-4-2\textmd{cosh}\big[\frac{1}{2}\sqrt{18{{\kappa}}^2-36 {\kappa}+12}c(x-ct)\big]\sqrt{-2{{\kappa}}^2+4{\kappa}}}.
\end{equation}

\begin{center}
\begin{tabular}{cc}
\epsfxsize=6.5cm \epsfysize=6.5cm \epsffile{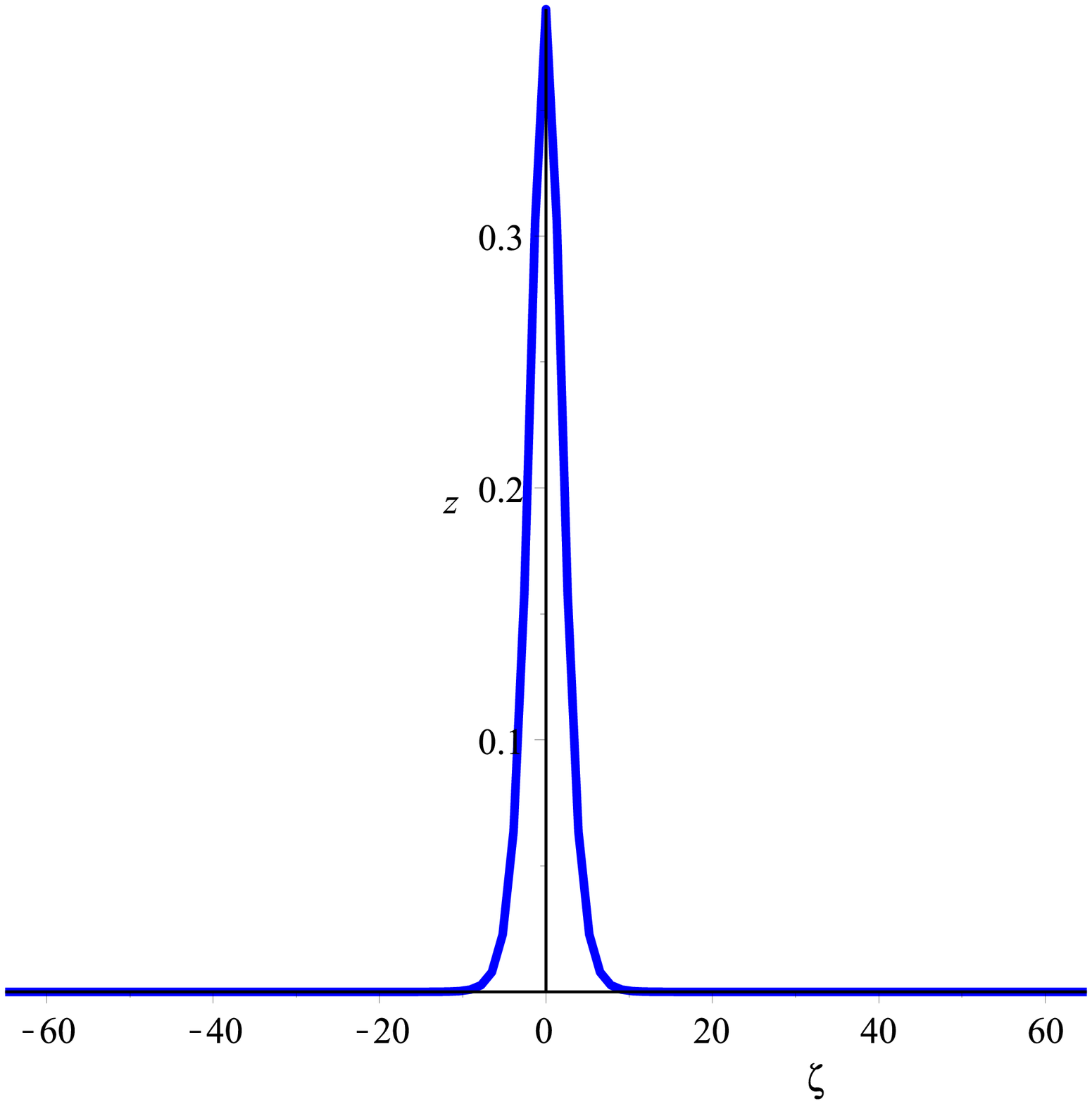}&
\epsfxsize=6.5cm \epsfysize=6.5cm \epsffile{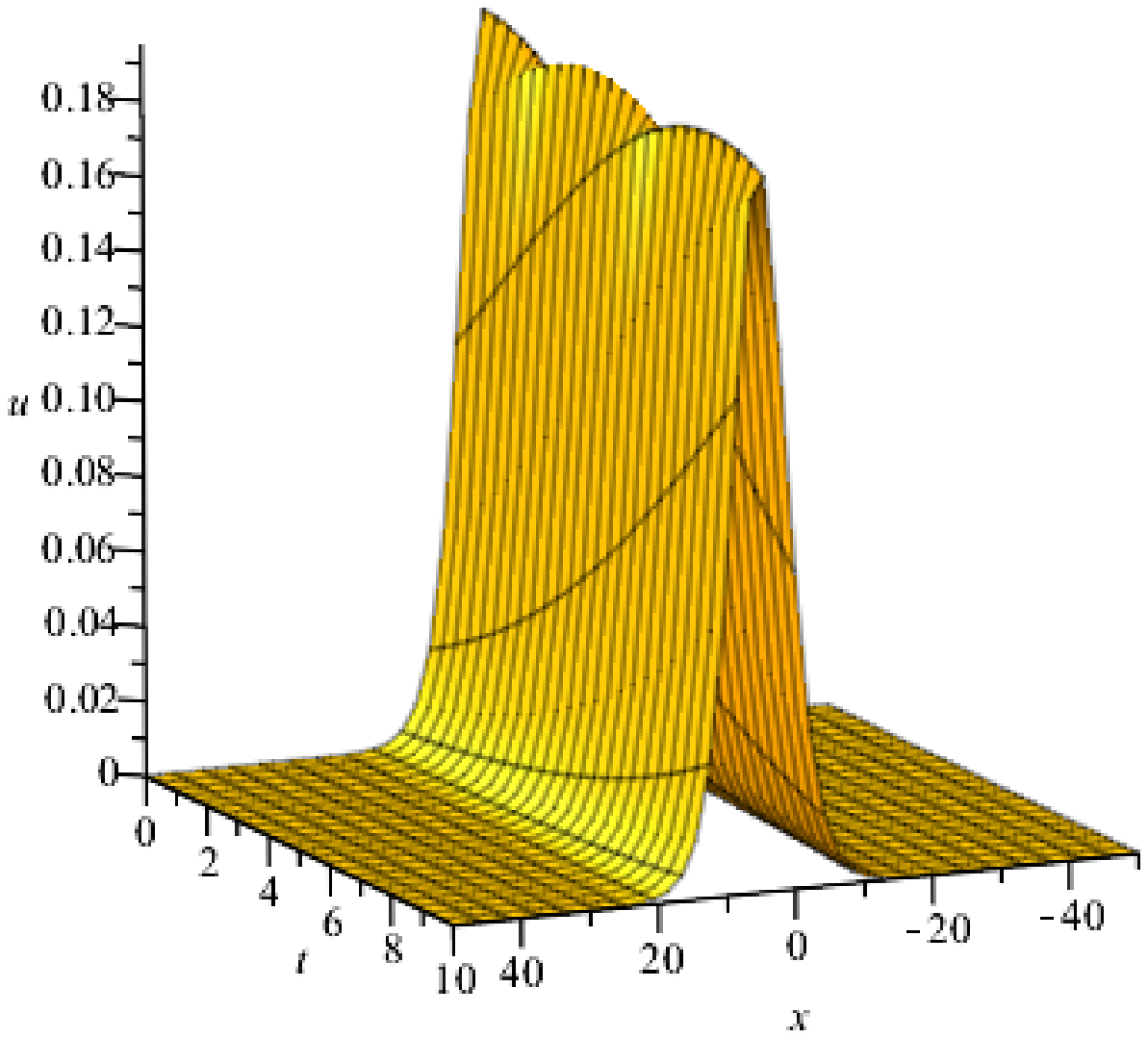}\\
\footnotesize{(a) Bright solitary wave } & \footnotesize{ (b) Exact solution }
\end{tabular}
\end{center}
\begin{center}
\footnotesize {{Fig. 3} \ \ Bright solitary wave of system $(\ref{ODE6})|_{\epsilon=0}$ and exact solution of Eq. $(\ref{wuzhang1})\mid_{\varepsilon=0}$ for $c=0.5$ and ${\kappa}=0.3$.}
\end{center}
In summary, we give the theorem as follows.
\begin{mythm}
For any wave speed $c>0$ and ${\kappa}\in (0,1-\frac{\sqrt3}{3})\cup (1+\frac{\sqrt3}{3},2)$, Eq. $(\ref{wuzhang1})\mid_{\varepsilon=0}$ has solitary wave solutions given by  (\ref{homoclinic2}) or (\ref{homoclinic4}).
\end{mythm}
\subsection{Kink and anti-kink wave solutions of Eq. $(\ref{wuzhang1})\mid_{\varepsilon=0}$}
According (\ref{Het}) and \textbf{Proposition} \ref{proposition} (\emph{iii}), we can obtain the expression of $y$ as follows:
\begin{equation}\label{het-expression1}
\begin{array}{rcl}
y=\pm \sqrt{\frac{3}{4}z^4-3z^3+3z^2}=\pm \sqrt{\frac{3}{4}(z-0)^2(2-z)^2}=\pm \frac{\sqrt3}{2}z(z-2).
\end{array}
\end{equation}
Similarly, due to $\frac{dz}{d\zeta}=y$, we have
\begin{equation}\label{het-expression2}
\zeta=\pm \int_{1}^{z}\frac{1}{\frac{\sqrt3}{2}s(s-2)}ds.
\end{equation}
For $\kappa=0$ or $\kappa=0$, the expression of the kink or anti-kink wave of system $(\ref{ODE6})|_{\epsilon=0}$ can be obtained by  (see Fig. 4(a) and (c))
\begin{equation}\label{homoclinic3}
z(\zeta)=-\frac{2c}{1+e^{\pm\sqrt3c\zeta}}.
\end{equation}
Therefore, we obtain the following kink and anti-kink wave solutions (see Fig. 4(b) and (b)):
\begin{equation}\label{het-expression3}
u(x,t)=-\frac{2c}{1+e^{\pm\sqrt3c(x-ct)}}.
\end{equation}

\begin{center}
\begin{tabular}{cc}
\epsfxsize=6.5cm \epsfysize=6.5cm \epsffile{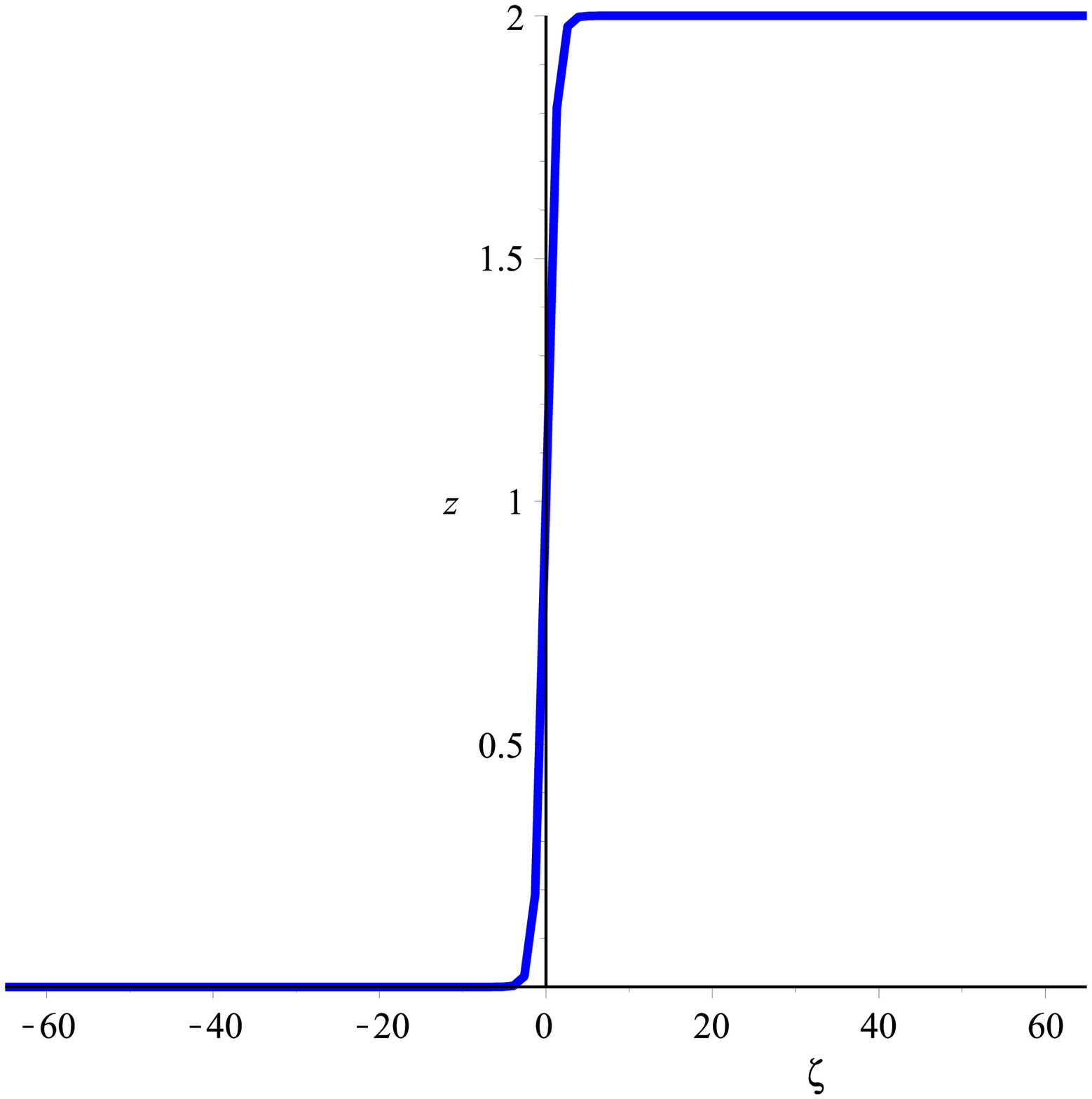}&
\epsfxsize=6.5cm \epsfysize=6.5cm \epsffile{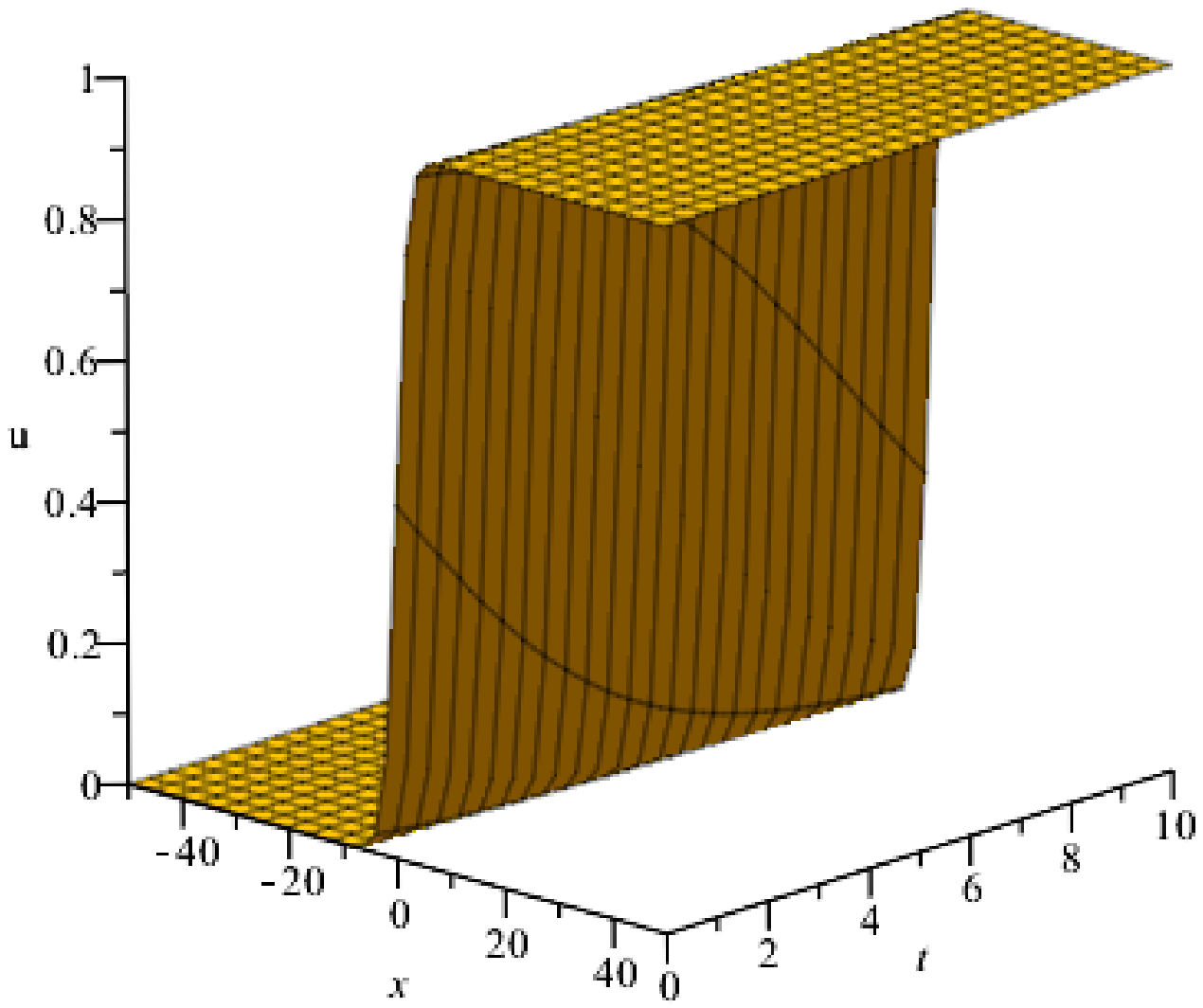}\\
\footnotesize{ (a) Kink wave} & \footnotesize{ (b) Kink wave solution}
\end{tabular}
\end{center}
\begin{center}
\begin{tabular}{cc}
\epsfxsize=6.5cm \epsfysize=6.5cm \epsffile{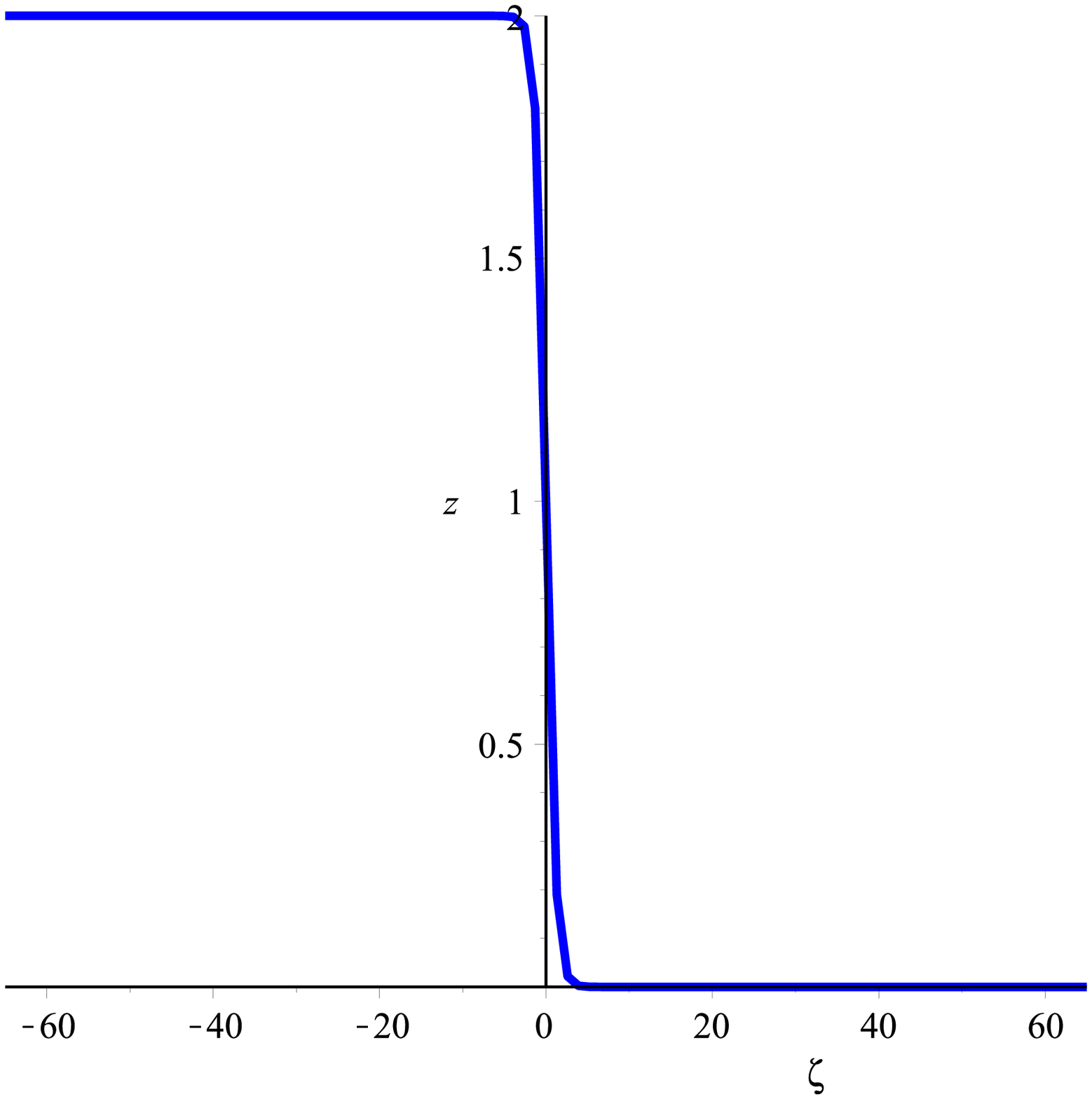}&
\epsfxsize=6.5cm \epsfysize=6.5cm \epsffile{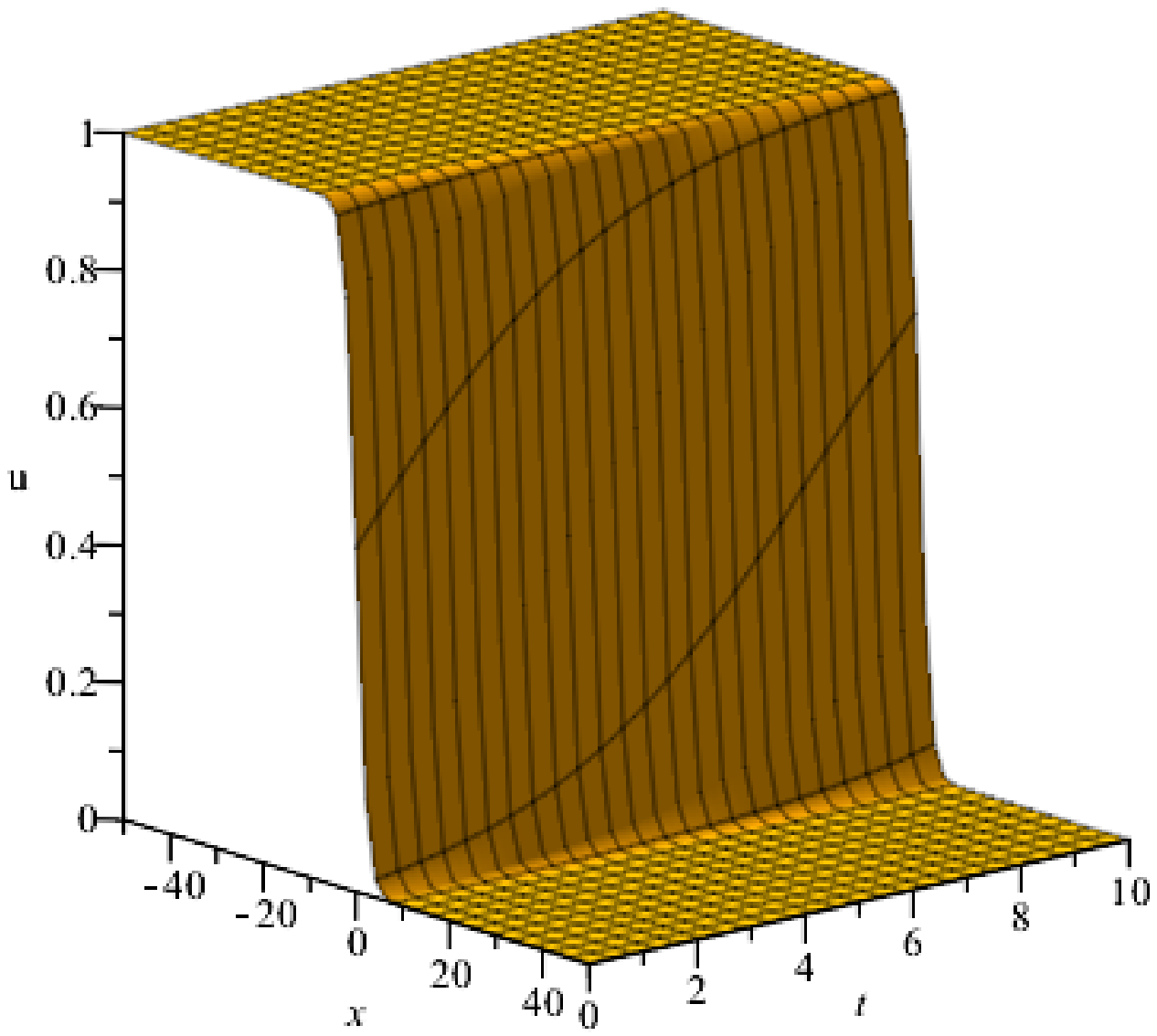}\\
\footnotesize{ (c) Anti-kink wave} & \footnotesize{ (d) Anti-kink wave solution}
\end{tabular}
\end{center}
\begin{center}
\footnotesize {{Fig. 4} \ \ Kink and anti-kink wave of system $(\ref{ODE6})|_{\epsilon=0}$. Kink and anti-kink wave solutions of Eq. $(\ref{wuzhang1})\mid_{\varepsilon=0}$ for $c=0.5$.}
\end{center}
Based on above analysis, we have the following theorems:
\begin{mythm}
For any wave speed $c>0$ and $\kappa=0$ or $\kappa=2$, Eq. $(\ref{wuzhang1})\mid_{\varepsilon=0}$ has kink or anti-kink solutions given by  (\ref{het-expression3}).
\end{mythm}
\section{ Bifurcations of the travelling wave solutions for the perturbed Eq. (\ref{wuzhang1})}
{In this section, before stating our main results on the existence of the solitary and kink wave solutions for the perturbed Eq. (\ref{wuzhang1}), denote
\begin{equation}\label{c-Melnikov1}
\begin{array}{rcl}
c_1(\kappa)\Big|_{\kappa\in(0,1-\frac{\sqrt3}{3})}&=&\frac{\sqrt{15}}{3}\big[\big(6\textmd{ln}2{{\kappa}}^3+6{{\kappa}}^3\textmd{ln}(\sqrt{6{{\kappa}}^2-12{\kappa}+4}+2{\kappa}-2)-18\textmd{ln}2{{\kappa}}^2\\
&&-18{{\kappa}}^2\textmd{ln}(\sqrt{6{{\kappa}}^2-12{\kappa}+4}+2{\kappa}-2)+12\textmd{ln}2{\kappa}\\
&&+12{\kappa}\textmd{ln}(\sqrt{6{{\kappa}}^2-12{\kappa}+4}+2{\kappa}-2)+2\sqrt{6{{\kappa}}^2-12{\kappa}+4}~ \big]\\
&&-\sqrt3 {\kappa}\big[2{{\kappa}}^2\textmd{ln}(-\sqrt{2{\kappa}-{{\kappa}}^2})+3\textmd{ln}2{{\kappa}}^2-6{\kappa}\textmd{ln}(-\sqrt{2{\kappa}-{{\kappa}}^2})-9\textmd{ln}2{\kappa}\\
&&+4\textmd{ln}(-\sqrt{2{\kappa}-{{\kappa}}^2})+6\textmd{ln}2\big)\big/\big(27{{\kappa}}^4 \sqrt{6{{\kappa}}^2-12{\kappa}+4}\\
&&+30 {{\kappa}}^3\textmd{ln}(-\sqrt{2{\kappa}-{{\kappa}}^2})+60{\kappa}\textmd{ln}(-\sqrt{2{\kappa}-{{\kappa}}^2})+15\textmd{ln}2{{\kappa}}^3\\
&&-30{{\kappa}}^3\textmd{ln}(\sqrt{6{{\kappa}}^2-12{\kappa}+4}+2{\kappa}-2)-90{{\kappa}}^2\textmd{ln}(-\sqrt{2{\kappa}-{{\kappa}}^2}) \\
&&+90{{\kappa}}^2\textmd{ln}(\sqrt{6{{\kappa}}^2-12{\kappa}+4}+2{\kappa}-2)+30\sqrt2{\kappa}\sqrt{(3{{\kappa}}^2-6{\kappa}+2)^3}\\
&&+294{{\kappa}}^2\sqrt{6{{\kappa}}^2-12{\kappa}+4}-45\textmd{ln}2{{\kappa}}^2-198{{\kappa}}^3\sqrt{6{{\kappa}}^2-12{\kappa}+4}\\
&&-60{\kappa}\textmd{ln}(\sqrt{6{{\kappa}}^2-12{\kappa}+4}+2{\kappa}-2)+10\sqrt2\sqrt{(3{{\kappa}}^2-6{\kappa}+2)^3}\\
&&-72{\kappa}\sqrt{6{{\kappa}}^2-12{\kappa}+4}+30\textmd{ln}2{\kappa}-18\sqrt{6{{\kappa}}^2-12{\kappa}+4}\ \big)\big]^{\frac{1}{2}}.
\end{array}
\end{equation}
and
\begin{equation}\label{c-Melnikov}
\begin{array}{rcl}
c_2(\kappa)\Big|_{\kappa\in(1+\frac{\sqrt3}{3},2)}&=&\frac{\sqrt{15}}{3}\big[\big(6\textmd{ln}2{{\kappa}}^3+6{{\kappa}}^3\textmd{ln}(\sqrt{6{{\kappa}}^2-12{\kappa}+4}+2{\kappa}-2)-18\textmd{ln}2{{\kappa}}^2\\
&&-18{{\kappa}}^2\textmd{ln}(\sqrt{6{{\kappa}}^2-12{\kappa}+4}+2{\kappa}-2)+12\textmd{ln}2{\kappa}\\
&&+12{\kappa}\textmd{ln}(\sqrt{6{{\kappa}}^2-12{\kappa}+4}+2{\kappa}-2)+2\sqrt{6{{\kappa}}^2-12{\kappa}+4}~ \big]\\
&&-\sqrt3 {\kappa}\big[{{\kappa}}^2\textmd{ln}(2{\kappa}-{{\kappa}}^2)+3\textmd{ln}2{{\kappa}}^2-3{\kappa}\textmd{ln}(2{\kappa}-{{\kappa}}^2)-9\textmd{ln}2{\kappa}\\
&&+2\textmd{ln}(2{\kappa}-{{\kappa}}^2)+6\textmd{ln}2\big)\big/\big(27{{\kappa}}^4 \sqrt{6{{\kappa}}^2-12{\kappa}+4}\\
&&+15 {{\kappa}}^3\textmd{ln}(2{\kappa}-{{\kappa}}^2)-198{{\kappa}}^3\sqrt{6{{\kappa}}^2-12{\kappa}+4}\\
&&-30{{\kappa}}^3\textmd{ln}(\sqrt{6{{\kappa}}^2-12{\kappa}+4}+2{\kappa}-2)+15\textmd{ln}2{{\kappa}}^3-45{{\kappa}}^2\textmd{ln}(2{\kappa}-{{\kappa}}^2) \\
&&+90{{\kappa}}^2\textmd{ln}(\sqrt{6{{\kappa}}^2-12{\kappa}+4}+2{\kappa}-2)+30\sqrt2{\kappa}\sqrt{(3{{\kappa}}^2-6{\kappa}+2)^3}\\
&&+294{{\kappa}}^2\sqrt{6{{\kappa}}^2-12{\kappa}+4}-45\textmd{ln}2{{\kappa}}^2+30{\kappa}\textmd{ln}(2{\kappa}-{{\kappa}}^2)\\
&&-60{\kappa}\textmd{ln}(\sqrt{6{{\kappa}}^2-12{\kappa}+4}+2{\kappa}-2)+10\sqrt2\sqrt{(3{{\kappa}}^2-6{\kappa}+2)^3}\\
&&-72{\kappa}\sqrt{6{{\kappa}}^2-12{\kappa}+4}+30\textmd{ln}2{\kappa}-18\sqrt{6{{\kappa}}^2-12{\kappa}+4}~ \big)\big]^{\frac{1}{2}},
\end{array}
\end{equation}

Now we are in a position to state our main results on the bifurcations of the travelling waves for the perturbed Eq. (\ref{wuzhang1}).

\begin{mythm}\label{theorem1}
(i) If ${\kappa}\in (0,1-\frac{\sqrt3}{3})$, then for sufficiently small $\epsilon$ ($0<\epsilon \ll 1$), there exists a wave speed $\tilde{c}_1(\kappa,\epsilon)=c_1({\kappa})+O(\epsilon)$ such that system (\ref{ODE6}) possesses a homoclinic orbit. Thus, Eq. (\ref{wuzhang1}) has a solitary wave solution $u=u(x,t,{\kappa},\epsilon)$ with a wave speed $c=\tilde{c}_1$.
\\

(ii) If ${\kappa}\in  (1+\frac{\sqrt3}{3},2)$, then for sufficiently small $\epsilon$ ($0<\epsilon \ll 1$), there exists a wave speed $\tilde{c}_2(\kappa,\epsilon)=c_2({\kappa})+O(\epsilon)$ such that system (\ref{ODE6}) possesses a homoclinic orbit. Thus, Eq. (\ref{wuzhang1}) has a solitary wave solution $u=u(x,t,{\kappa},\epsilon)$ with a wave speed $c=\tilde{c}_2$.
\\

 (iii) If $\kappa=0$ or $\kappa=2$, then for sufficiently small $\epsilon$ ($0<\epsilon \ll 1$), there exists a unique wave speed $\tilde{c}_3(\epsilon)=\sqrt{15}/3+O(\epsilon)$ such that system
(\ref{ODE6}) has a heteroclinic orbit. Therefore,  Eq. (\ref{wuzhang1}) has a pair of kink and anti-kink wave solutions $u=u(x,t,\kappa,\epsilon)$ with a wave speed $c=\tilde{c}_3$.
\end{mythm}

\noindent
{\bf Remark 1.}
Different kinds of travelling wave solutions exist in the different bifurcation parameter regions. For DWLE, solitary wave and kink wave can not coexist at a same wave speed.\\

\noindent
{\bf Remark 2.}
The exact limit wave speed is given for each bifurcation parameter regions.\\
}

In order to prove \textbf{Theorem \ref{theorem1}}, we need to introduce some lemmas.

Poincar\'{e} map $d_i$ ($i=1,2$) is a powerful tool to study the bifurcation of homoclinic or heteroclinic orbits (e.g., see \cite{Hanmaoan-book,WS-NY,ZK-ND,wen-MMAS}). \\
\textbf{Case a}. It is defined by (see Fig. 5(a))
$$d_1:A(h)\longrightarrow P(h,c,\varepsilon),$$
and
\begin{equation*}\label{Melnikov}
\begin{array}{rcl}
d_1(h,c,\varepsilon)&=&\displaystyle \int_{A(h)}^{P(h,c,\varepsilon)}dH\\
&=&\displaystyle \int_{-\infty}^{+\infty}(\frac{\partial H}{\partial x}\frac{dx}{dt}+\frac{\partial H}{\partial y}\frac{dy}{dt})\mid_{L_{(h,c,\varepsilon)}}dt,
\end{array}
\end{equation*}
where $A(h)$ is the initial point and $P(h,c,\varepsilon)$ is the mapping point. $L(h,c,\varepsilon)$ and $L(h,c)$ are the perturbed and unperturbed orbits, respectively. We have $\lim\limits_{\varepsilon\rightarrow 0 }L(h,c,\varepsilon)=L(h,c)$. And $d_1(h,c,\varepsilon)$ can be Taylor expanded in $\varepsilon$, one has
$$d_1(h,c,\varepsilon)=\varepsilon M_{\textmd{hom}}(h,c)+O(\varepsilon^2).$$
{\textbf{Case b}. A heteroclinic orbit $L(h,c)$ connecting two saddle points $S_1$ and $S_2$. Then, $L^{+}_{s}(h,c,\varepsilon)$ and $L^{-}_{u}(h,c,\varepsilon)$ represent the stable and unstable manifolds of perturbed heteroclinic orbit $L(h,c,\varepsilon)$. Taking a point $P(h,c)\in L(h,c)$, we let $L^*$ be a segment normal of $L(h,c)$ at point $P(h,c)$. For $0<\varepsilon\ll1$, we assume that $L^{+}_{s}(h,c,\varepsilon)$ and $L^{-}_{u}(h,c,\varepsilon)$ intersect the normal line $L^*$  transversally at points $P^{+}_{s}(h,c,\varepsilon)$ and $P^{-}_{u}(h,c,\varepsilon)$ (see Fig. 5(b)).

Let
\begin{equation*}\label{Melnikov}
\begin{array}{rcl}
d_2(h,c,\varepsilon)&=&-\overrightarrow{n}\cdot\overrightarrow{P^{+}_{s}P^{-}_{u}},
\end{array}
\end{equation*}
where $\overrightarrow{n}=\frac{(H_z(P(h,c)),H_y(P(h,c)))}{\mid(H_y(P(h,c)),-H_z(P(h,c)))\mid}$.  And $d_2(h,c,\varepsilon)$ can be Taylor expanded in $\varepsilon$, we have
$$d_2(h,c,\varepsilon)=\varepsilon \cdot A\cdot M_{\textmd{het}}(h,c)+O(\varepsilon^2),$$
where $A$ is a constant.

Usually, $d_i(h,c,\varepsilon)$ ($i=1,2$) is used to measure the distance between perturbed stable and unstable manifolds. Then, we say that $M_{\textmd{hom}}(h,c)$ or $M_{\textmd{het}}(h,c)$ is the so called ``first-order" Melnikov integral.
\begin{center}
\begin{tabular}{cc}
\epsfxsize=5.5cm \epsfysize=3.5cm \epsffile{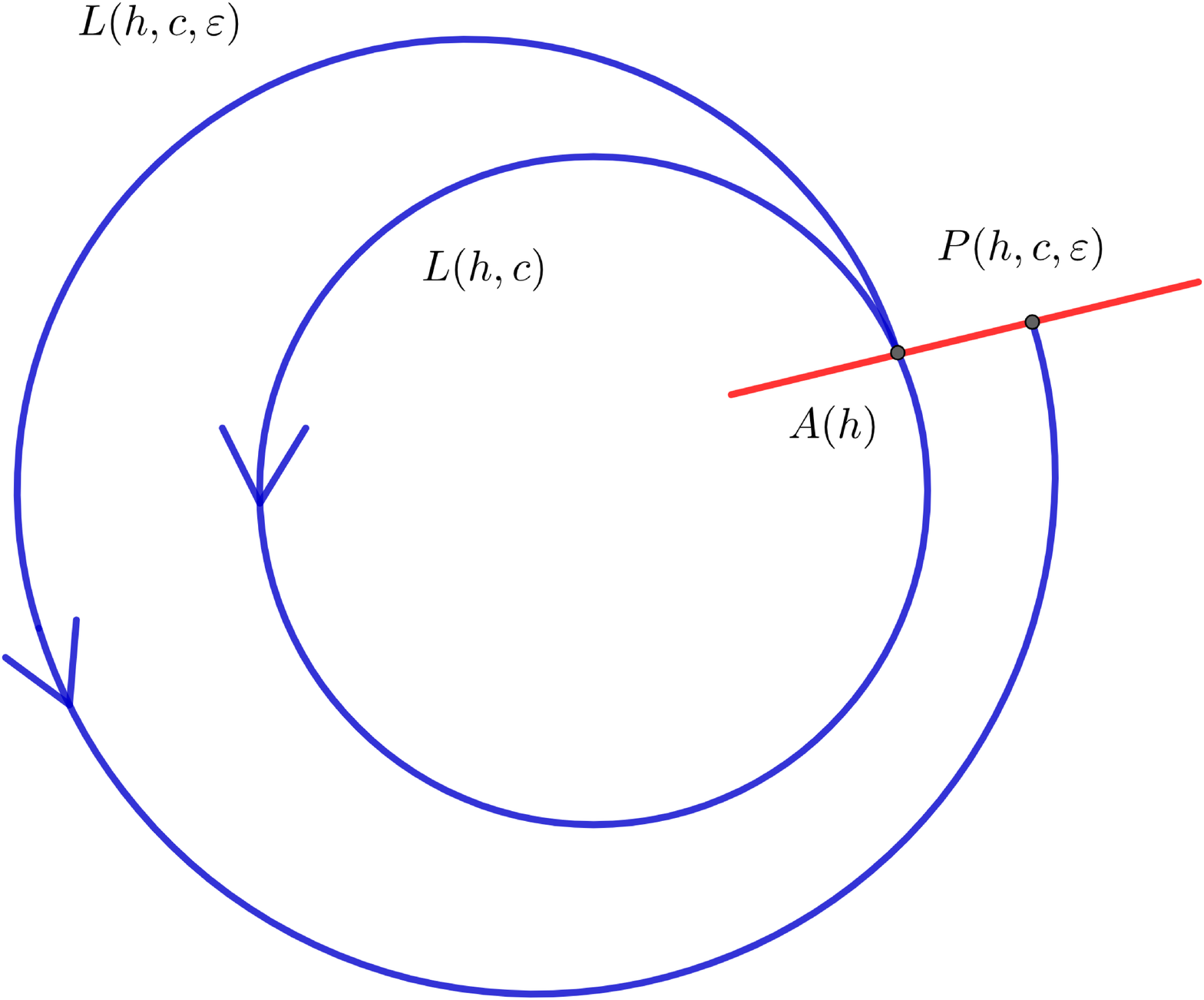}&
\epsfxsize=6.5cm \epsfysize=4cm \epsffile{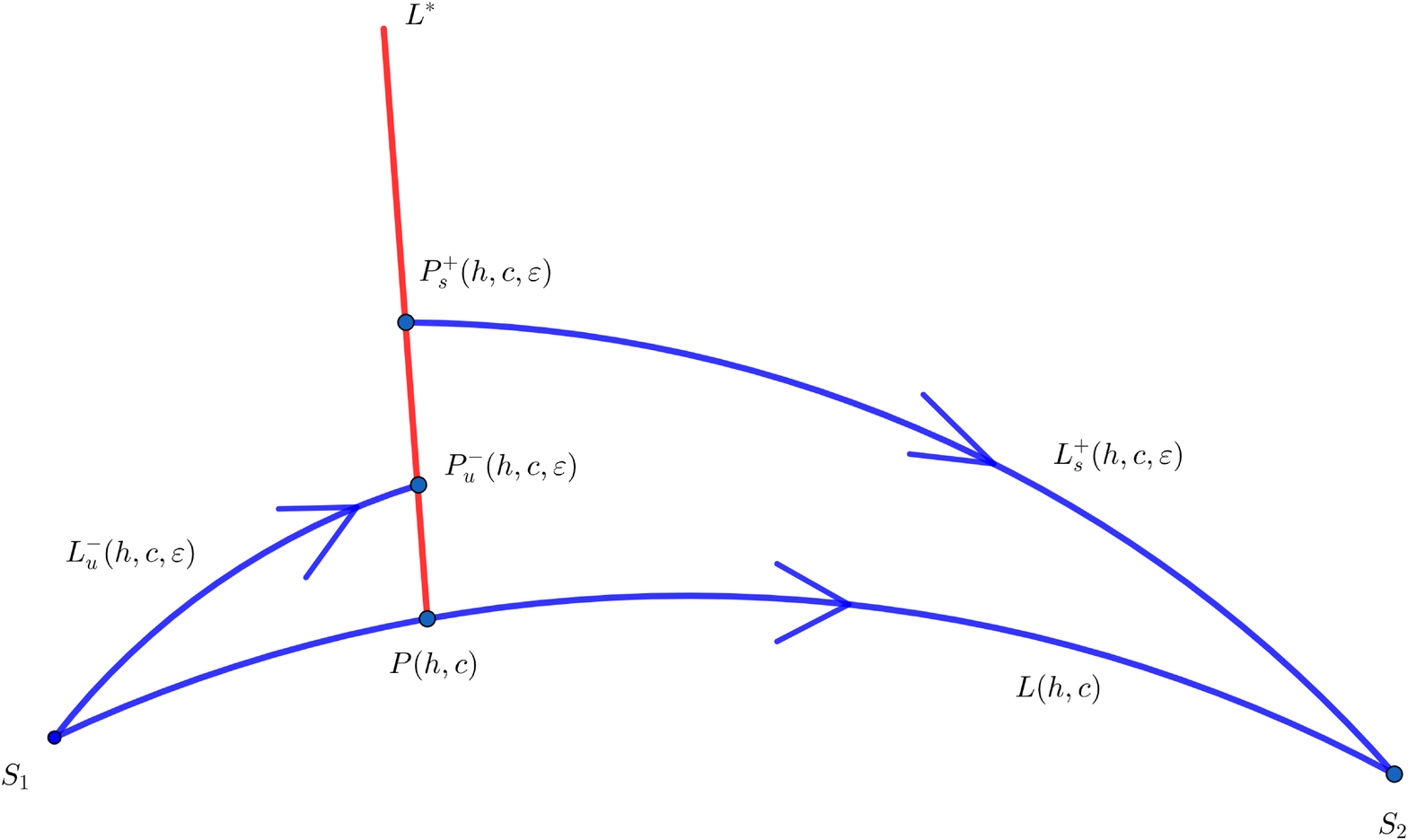}\\
\footnotesize{ (a) The type of periodic or homoclinic orbits } & \footnotesize{ (b) The type of heteroclinic orbits}
\end{tabular}
\end{center}
\begin{center}
\footnotesize {{Fig. 5} \ \  Poincar\'{e} map }
\end{center}
In this paper, ``first-order" Melnikov integral is employed to detect the persistence of the homoclinic and heteroclinic orbits under small perturbation, respectively.
Consequently, the existence of solitary and kink wave solutions are proved under small perturbation.}

Firstly, the homoclinic Melnikov function of system (\ref{ODE6}) is defined by
\begin{equation}\label{Melnikov}
\begin{array}{rcl}
M_{\textmd{hom}}(c,{\kappa})=\displaystyle\oint_{\Gamma({\kappa})}(\frac{9}{2}z^2-9z+3+\frac{1}{c^2})y^2d\zeta=\frac{1}{c^2}I({\kappa})+J({\kappa}),
\end{array}
\end{equation}
where $I({\kappa})=\displaystyle\oint_{\Gamma({\kappa})}y^2d\zeta$ and $J({\kappa})=\displaystyle\oint_{\Gamma({\kappa})}(\frac{9}{2}z^2-9z+3)y^2d\zeta$.
\begin{mylem}\label{lemma1}
For any ${\kappa}\in (0,1-\frac{\sqrt3}{3})\cup (1+\frac{\sqrt3}{3},2)$, there exists a positive root $c=c({\kappa})$ of $M_{\textmd{hom}}(c,{\kappa})=0$. Moreover, $\frac{\partial M_{\textmd{hom}}}{\partial c}\mid_{c=c({\kappa})}\neq0$.
\end{mylem}
\emph{Proof. } Obviously, we know that $I({\kappa})=\displaystyle\oint_{\Gamma({\kappa})}y^2d\zeta>0$. Because $J({\kappa})$ integrates along the homoclinic orbit $\Gamma({\kappa})$, then
\begin{equation}\label{integrate}
\begin{array}{rcl}
y''=(\frac{9}{2}z^2-9z+3)z'=(\frac{9}{2}z^2-9z+3)y,
\end{array}
\end{equation}
which yields
\begin{equation}\label{integrate1}
\begin{array}{rcl}
dy'=(\frac{9}{2}z^2-9z+3)z'd\zeta=(\frac{9}{2}z^2-9z+3)yd\zeta.
\end{array}
\end{equation}
By integration of parts, it is not difficult to have
\begin{equation}\label{integrate2}
\begin{array}{rcl}
J({\kappa})&=&\displaystyle\oint_{\Gamma({\kappa})}(\frac{9}{2}z^2-9z+3)y^2d\zeta=\displaystyle\oint_{\Gamma({\kappa})}ydy'\\
&=&-\displaystyle\int_\mathbb{R}(y')^2d\zeta<0.
\end{array}
\end{equation}
Hence, $-\frac{I({\kappa})}{J({\kappa})}>0$, there exists a single positive root $c=c({\kappa})$ of $M_{\textmd{hom}}(c,{\kappa})=0$ that
\begin{equation}\label{integrate3}
\begin{array}{rcl}
c({\kappa})=\sqrt{-\frac{I({\kappa})}{J({\kappa})}}.
\end{array}
\end{equation}
On the other hand, it follows from (\ref{Melnikov}) that
\begin{equation}\label{integrate4}
\begin{array}{rcl}
\frac{\partial M_{\textmd{hom}}(c,{\kappa})}{\partial c}\mid_{c=c({\kappa})}=-2\frac{1}{c^3}I({\kappa})=-2J({\kappa})\sqrt{-\frac{J({\kappa})}{I({\kappa})}}>0.
\end{array}
\end{equation}
The proof is completed.

In fact, the analytical expressions of $I({\kappa})$ and $J({\kappa})$ can be computed by dividing them into two cases.

\textbf{Case (I)}: $1+\frac{\sqrt3}{3}<{\kappa}<2$.

Firstly, since $I{({\kappa})}$ and $J{({\kappa})}$ are integrated over a closed curve, and the time variable $\zeta$ can be represented by the state variable $z$ on homoclinic orbit. Thus, by (\ref{expression1}), we have
\begin{equation}\label{I}
\begin{array}{rcl}
I({\kappa})&=&\displaystyle\oint_{\Gamma({\kappa})}y^2d\zeta=\displaystyle\oint_{\Gamma({\kappa})}ydz\\
&=&2 \displaystyle\int_{z_+}^{{\kappa}}\sqrt{\frac{3}{4}z^4-3z^3+3z^2-(3{{\kappa}}^3-9{{\kappa}}^2+6{\kappa})z+\frac{9}{4}{{\kappa}}^4-6{{\kappa}}^3+3{{\kappa}}^2}dz\\
&=&2 \displaystyle\int_{z_+}^{{\kappa}}\sqrt{\frac{3}{4}(z-{\kappa})^2(z-z_+)(z-z_-)}dz\\
&=&\sqrt{3} \big[\displaystyle\int_{z_+}^{{\kappa}}{\kappa}\sqrt{(z-z_+)(z-z_-)}dz-\displaystyle\int_{z_+}^{{\kappa}}z\sqrt{(z-z_+)(z-z_-)}dz \big]\\
&=& \sqrt{3} \big [{\kappa}I_1({\kappa})-I_2({\kappa}) \big],
\end{array}
\end{equation}
where
$I_1({\kappa})=\displaystyle\int_{z_+}^{{\kappa}}\sqrt{(z-z_+)(z-z_-)}dz$ and $I_2({\kappa})=\displaystyle\int_{z_+}^{{\kappa}}z\sqrt{(z-z_+)(z-z_-)}dz$.

 Note that the formulas 2.261, 2.262 (1) and 2.262 (2) in handbook \cite{GI-SV} read:
\begin{equation}\label{I1}
\begin{array}{rcl}
\displaystyle\int \frac{dx}{\sqrt R}&=& \frac{1}{\sqrt{c_1}}\textmd{ln}(2\sqrt{c_1R}+2c_1x+b), \ \ (c_1>0, \ 2c_1x+b>\sqrt{-\Delta}, \ \Delta<0),\\
\displaystyle\int\sqrt R dx&=& \frac{(2c_1x+b)\sqrt R}{4c_1}+\frac{\Delta}{8c_1}\displaystyle \int \frac{dx}{\sqrt R},\\
\displaystyle\int x\sqrt R dx&=& \frac{\sqrt{R^3}}{3c_1}-\frac{(2c_1x+b)b\sqrt R}{8c_1^2}-\frac{b\Delta}{16c_1^2}\displaystyle \int \frac{dx}{\sqrt R},
\end{array}
\end{equation}
where $R=a+bx+c_1x^2$ and $\Delta=4ac_1-b^2$.
Combining (\ref{I}) and (\ref{I1}), it follows that
\begin{equation}\label{I2}
\begin{array}{rcl}
I({\kappa})&=&\frac{\sqrt3}{3}\big[6\textmd{ln}2{{\kappa}}^3+6{{\kappa}}^3\textmd{ln}(\sqrt{6{{\kappa}}^2-12{\kappa}+4}+2{\kappa}-2)-18\textmd{ln}2{{\kappa}}^2\\
&&-18{{\kappa}}^2\textmd{ln}(\sqrt{6{{\kappa}}^2-12{\kappa}+4}+2{\kappa}-2)+12\textmd{ln}2{\kappa}\\
&&+12{\kappa}\textmd{ln}(\sqrt{6{{\kappa}}^2-12{\kappa}+4}+2{\kappa}-2)+2\sqrt{6{{\kappa}}^2-12{\kappa}+4}~ \big]\\
&&-\sqrt3 {\kappa}\big[{{\kappa}}^2\textmd{ln}(2{\kappa}-{{\kappa}}^2)+3\textmd{ln}2{{\kappa}}^2-3{\kappa}\textmd{ln}(2{\kappa}-{{\kappa}}^2)-9\textmd{ln}2{\kappa}\\
&&+2\textmd{ln}(2{\kappa}-{{\kappa}}^2)+6\textmd{ln}2\big].
\end{array}
\end{equation}
Secondly,
\begin{equation}\label{J}
\begin{array}{rcl}
J({\kappa})&=&\displaystyle\oint_{\Gamma({\kappa})}(\frac{9}{2}z^2-9z+3)y^2d\zeta
=\displaystyle\oint_{\Gamma({\kappa})}(\frac{9}{2}z^2-9z+3)ydz\\
&=&2 \displaystyle\int_{z_+}^{{\kappa}}(\frac{9}{2}z^2-9z+3)\scriptstyle\sqrt{\frac{3}{4}z^4-3z^3+3z^2-(3{{\kappa}}^3-9{{\kappa}}^2+6{\kappa})z+\frac{9}{4}{{\kappa}}^4-6{{\kappa}}^3+3{{\kappa}}^2}\displaystyle dz\\
&=&\sqrt{3} \displaystyle\int_{z_+}^{{\kappa}}(\frac{9}{2}z^2-9z+3)\sqrt{({\kappa}-z)^2(z-z_+)(z-z_-)}dz\\
&=&\sqrt{3}\displaystyle\int_{z_+}^{{\kappa}}\big[-\frac{9}{2}z^3+(\frac{9}{2}{\kappa}+9)z^2-(9{\kappa}+3)z+3{\kappa}\big]\sqrt{(z-z_+)(z-z_-)}dz \\
&=& \sqrt{3} \big [-\frac{9}{2}J_1({\kappa})+(\frac{9}{2}{\kappa}+9)J_2({\kappa})-(9{\kappa}+3)I_2({\kappa})+3{\kappa}I_1({\kappa}) \big],
\end{array}
\end{equation}
where
$J_1({\kappa})=\displaystyle\int_{z_+}^{{\kappa}}z^3\sqrt{(z-z_+)(z-z_-)}dz$, $J_2({\kappa})=\displaystyle\int_{z_+}^{{\kappa}}z^2\sqrt{(z-z_+)(z-z_-)}dz$.

Similarly, by the formulas 2.262 (3) and 2.262 (4) in \cite{GI-SV}, they are expressed by:
\begin{equation}\label{J1}
\begin{array}{rcl}
\displaystyle\int x^2\sqrt R dx&=& (\frac{x}{4c_1}-\frac{5b}{24c_1^2})\sqrt {R^3}+(\frac{5b^2}{16c_1^2}-\frac{a}{4c_1})\frac{(2c_1x+b)\sqrt R}{4c_1}+(\frac{5b^2}{16c_1^2}-\frac{a}{4c_1})\frac{\Delta}{8c_1}\displaystyle \int \frac{dx}{\sqrt R},\\
\displaystyle\int x^3\sqrt R dx&=& (\frac{x^2}{5c_1}-\frac{7bx}{40c_1^2}+\frac{7b^2}{48c_1^3}-\frac{2a}{15c_1^2})\sqrt {R^3}
-(\frac{7b^3}{32c_1^3}-\frac{3ab}{8c_1^2})\frac{(2c_1x+b)\sqrt R}{4c_1}\\
&&- (\frac{7b^3}{32c_1^3}-\frac{3ab}{8c_1^2})\frac{\Delta}{8c_1}\displaystyle \int \frac{dx}{\sqrt R}.
\end{array}
\end{equation}
Using (\ref{I1}), (\ref{J}) and (\ref{J1}), we have
\begin{equation}\label{J2}
\begin{array}{rcl}
J({\kappa})&=&-\frac{\sqrt3}{5} \big[27{{\kappa}}^4 \sqrt{6{{\kappa}}^2-12{\kappa}+4}+15 {{\kappa}}^3\textmd{ln}(2{\kappa}-{{\kappa}}^2)-198{{\kappa}}^3\sqrt{6{{\kappa}}^2-12{\kappa}+4}\\
&&-30{{\kappa}}^3\textmd{ln}(\sqrt{6{{\kappa}}^2-12{\kappa}+4}+2{\kappa}-2)+15\textmd{ln}2{{\kappa}}^3-45{{\kappa}}^2\textmd{ln}(2{\kappa}-{{\kappa}}^2) \\
&&+90{{\kappa}}^2\textmd{ln}(\sqrt{6{{\kappa}}^2-12{\kappa}+4}+2{\kappa}-2)+30\sqrt2{\kappa}\sqrt{(3{{\kappa}}^2-6{\kappa}+2)^3}\\
&&+294{{\kappa}}^2\sqrt{6{{\kappa}}^2-12{\kappa}+4}-45\textmd{ln}2{{\kappa}}^2+30{\kappa}\textmd{ln}(2{\kappa}-{{\kappa}}^2)\\
&&-60{\kappa}\textmd{ln}(\sqrt{6{{\kappa}}^2-12{\kappa}+4}+2{\kappa}-2)+10\sqrt2\sqrt{(3{{\kappa}}^2-6{\kappa}+2)^3}\\
&&-72{\kappa}\sqrt{6{{\kappa}}^2-12{\kappa}+4}+30\textmd{ln}2{\kappa}-18\sqrt{6{{\kappa}}^2-12{\kappa}+4}~ \big].
\end{array}
\end{equation}

\textbf{Case (II)}: $0<{\kappa}<1-\frac{\sqrt3}{3}$.

Same calculation method as case (I), we obtain
\begin{equation}\label{I22}
\begin{array}{rcl}
I({\kappa})&=&\displaystyle\oint_{\Gamma({\kappa})}y^2d\zeta=\displaystyle\oint_{\Gamma({\kappa})}ydz=2 \displaystyle\int_{z_-}^{{\kappa}}\cdots dz\\
&=&\frac{\sqrt3}{3}\big[6\textmd{ln}2{{\kappa}}^3+6{{\kappa}}^3\textmd{ln}(\sqrt{6{{\kappa}}^2-12{\kappa}+4}+2{\kappa}-2)-18\textmd{ln}2{{\kappa}}^2\\
&&-18{{\kappa}}^2\textmd{ln}(\sqrt{6{{\kappa}}^2-12{\kappa}+4}+2{\kappa}-2)+12\textmd{ln}2{\kappa}\\
&&+12{\kappa}\textmd{ln}(\sqrt{6{{\kappa}}^2-12{\kappa}+4}+2{\kappa}-2)+2\sqrt{6{{\kappa}}^2-12{\kappa}+4}~ \big]\\
&&-\sqrt3 {\kappa}\big[2{{\kappa}}^2\textmd{ln}(-\sqrt{2{\kappa}-{{\kappa}}^2})+3\textmd{ln}2{{\kappa}}^2-6{\kappa}\textmd{ln}(-\sqrt{2{\kappa}-{{\kappa}}^2})-9\textmd{ln}2{\kappa}\\
&&+4\textmd{ln}(-\sqrt{2{\kappa}-{{\kappa}}^2})+6\textmd{ln}2\big],
\end{array}
\end{equation}
and
\begin{equation}\label{J22}
\begin{array}{rcl}
J({\kappa})&=&\displaystyle\oint_{\Gamma({\kappa})}(\frac{9}{2}z^2-9z+3)y^2d\zeta=\displaystyle\oint_{\Gamma({\kappa})}(\frac{9}{2}z^2-9z+3)ydz=2 \displaystyle\int_{z_-}^{{\kappa}}\cdots dz\\
&=&-\frac{\sqrt3}{5} \big[27{{\kappa}}^4 \sqrt{6{{\kappa}}^2-12{\kappa}+4}+30 {{\kappa}}^3\textmd{ln}(-\sqrt{2{\kappa}-{{\kappa}}^2})+60{\kappa}\textmd{ln}(-\sqrt{2{\kappa}-{{\kappa}}^2})\\
&&-30{{\kappa}}^3\textmd{ln}(\sqrt{6{{\kappa}}^2-12{\kappa}+4}+2{\kappa}-2)+15\textmd{ln}2{{\kappa}}^3-90{{\kappa}}^2\textmd{ln}(-\sqrt{2{\kappa}-{{\kappa}}^2}) \\
&&+90{{\kappa}}^2\textmd{ln}(\sqrt{6{{\kappa}}^2-12{\kappa}+4}+2{\kappa}-2)+30\sqrt2{\kappa}\sqrt{(3{{\kappa}}^2-6{\kappa}+2)^3}\\
&&+294{{\kappa}}^2\sqrt{6{{\kappa}}^2-12{\kappa}+4}-45\textmd{ln}2{{\kappa}}^2-198{{\kappa}}^3\sqrt{6{{\kappa}}^2-12{\kappa}+4}\\
&&-60{\kappa}\textmd{ln}(\sqrt{6{{\kappa}}^2-12{\kappa}+4}+2{\kappa}-2)+10\sqrt2\sqrt{(3{{\kappa}}^2-6{\kappa}+2)^3}\\
&&-72{\kappa}\sqrt{6{{\kappa}}^2-12{\kappa}+4}+30\textmd{ln}2{\kappa}-18\sqrt{6{{\kappa}}^2-12{\kappa}+4}\ \big].
\end{array}
\end{equation}

We plot the algebraic curve of $I({\kappa})$ with respect to ${\kappa}$ and $J({\kappa})$ with respect to ${\kappa}$, respectively. (see Fig. 6 and Fig. 7)
\begin{center}
\begin{tabular}{cc}
\epsfxsize=5.5cm \epsfysize=5.5cm \epsffile{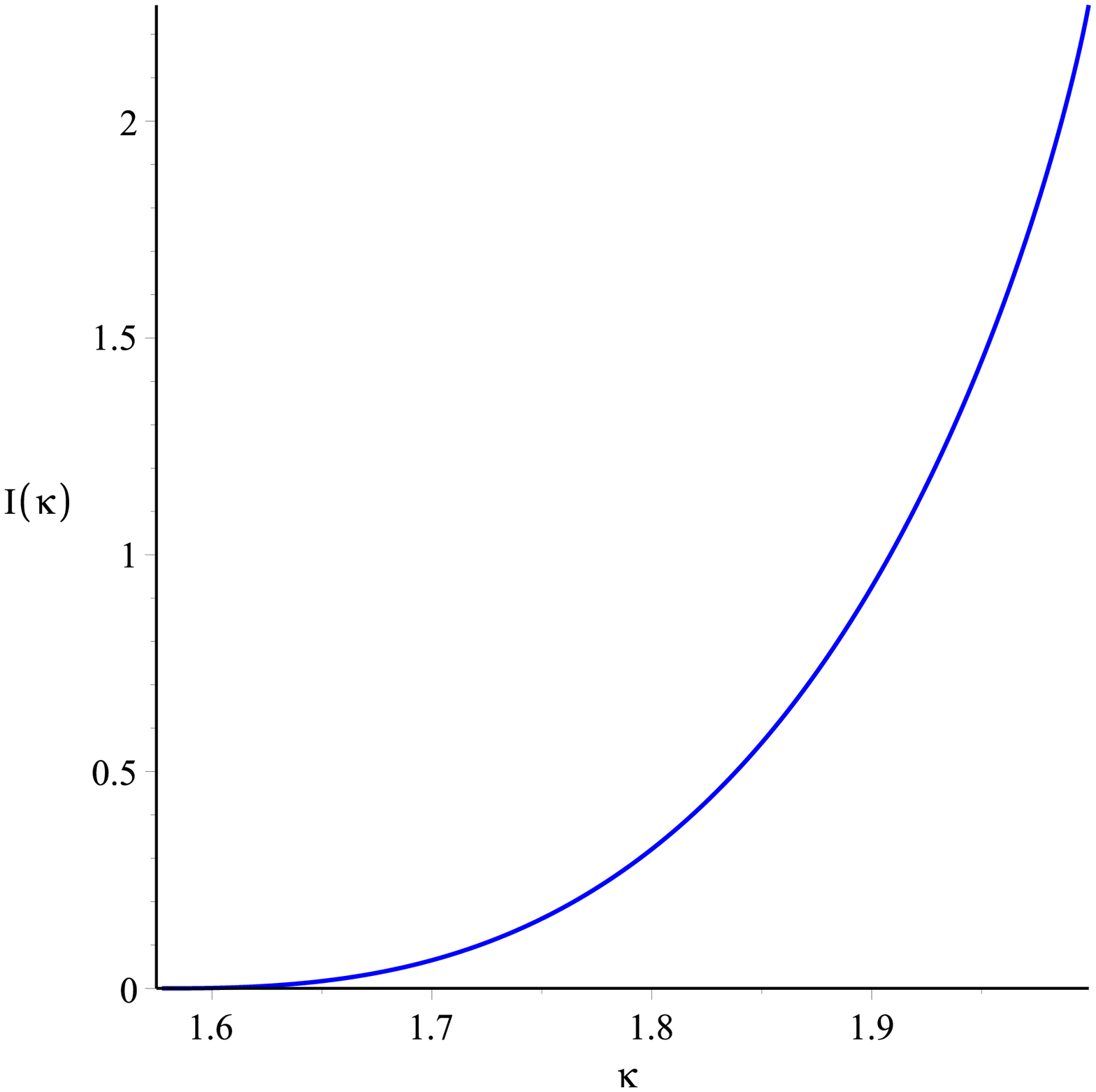}&
\epsfxsize=5.5cm \epsfysize=5.5cm \epsffile{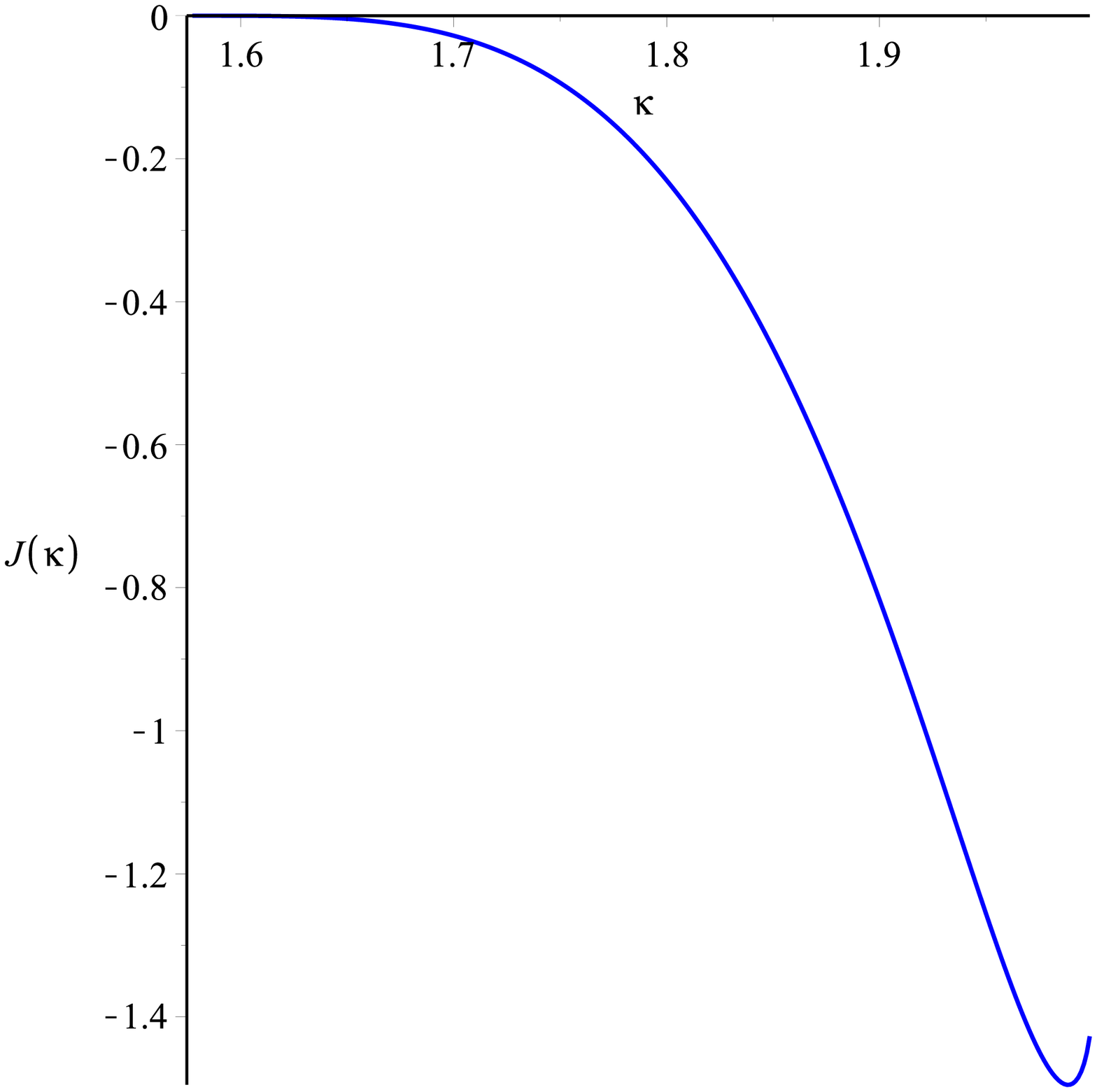}\\
\footnotesize{ (a) (${\kappa}$,\ $I(\kappa)$)} & \footnotesize{ (b) (${\kappa}$,\ $J(\kappa)$)}
\end{tabular}
\end{center}
\begin{center}
\footnotesize {{Fig. 6} \ \  The algebraic curves of $I({\kappa})$ and $J({\kappa})$ with respect to ${\kappa}$ for $1+\frac{\sqrt3}{3}<{\kappa}<2$.}
\end{center}
\begin{center}
\begin{tabular}{cc}
\epsfxsize=5.5cm \epsfysize=5.5cm \epsffile{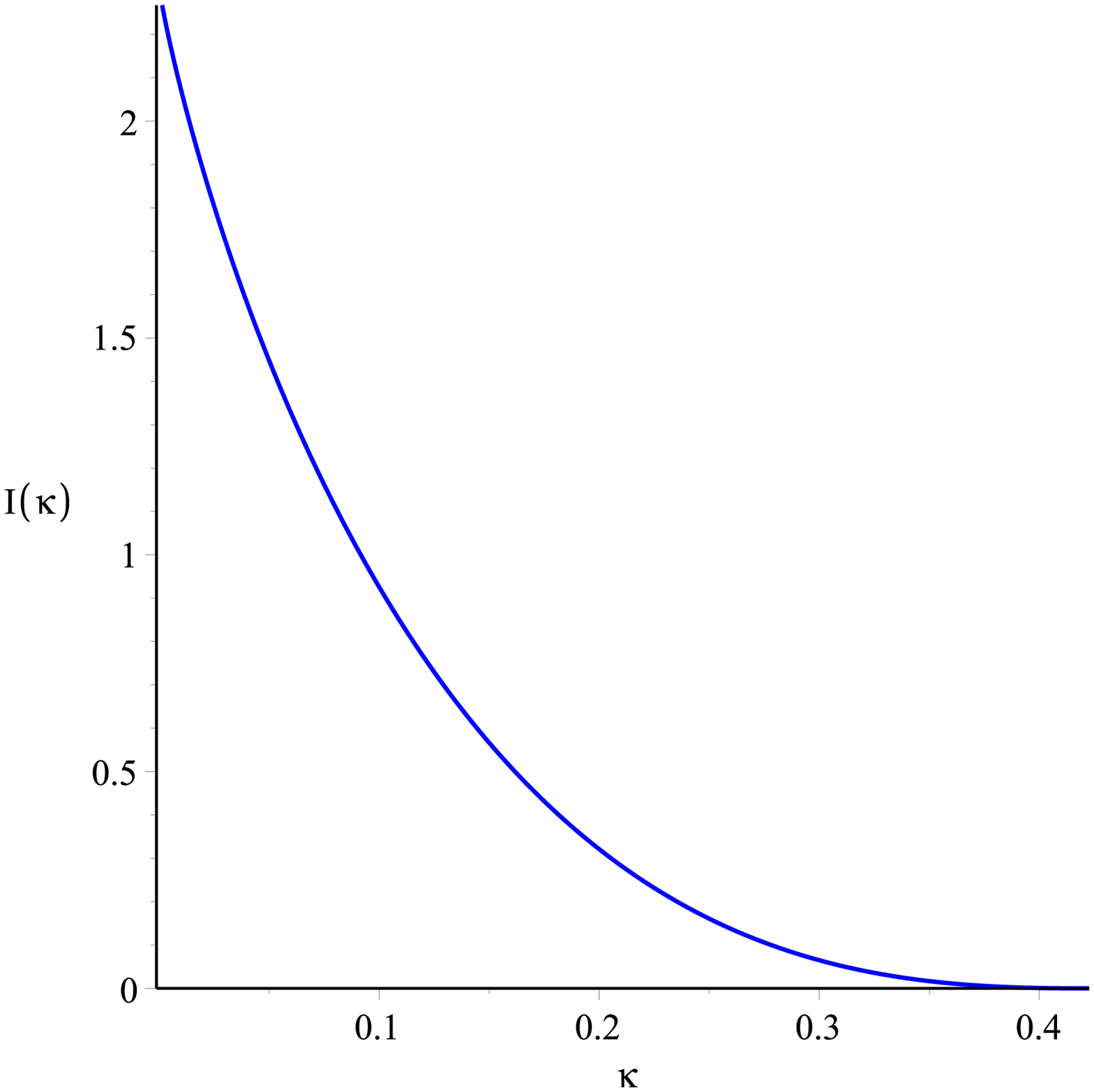}&
\epsfxsize=5.5cm \epsfysize=5.5cm \epsffile{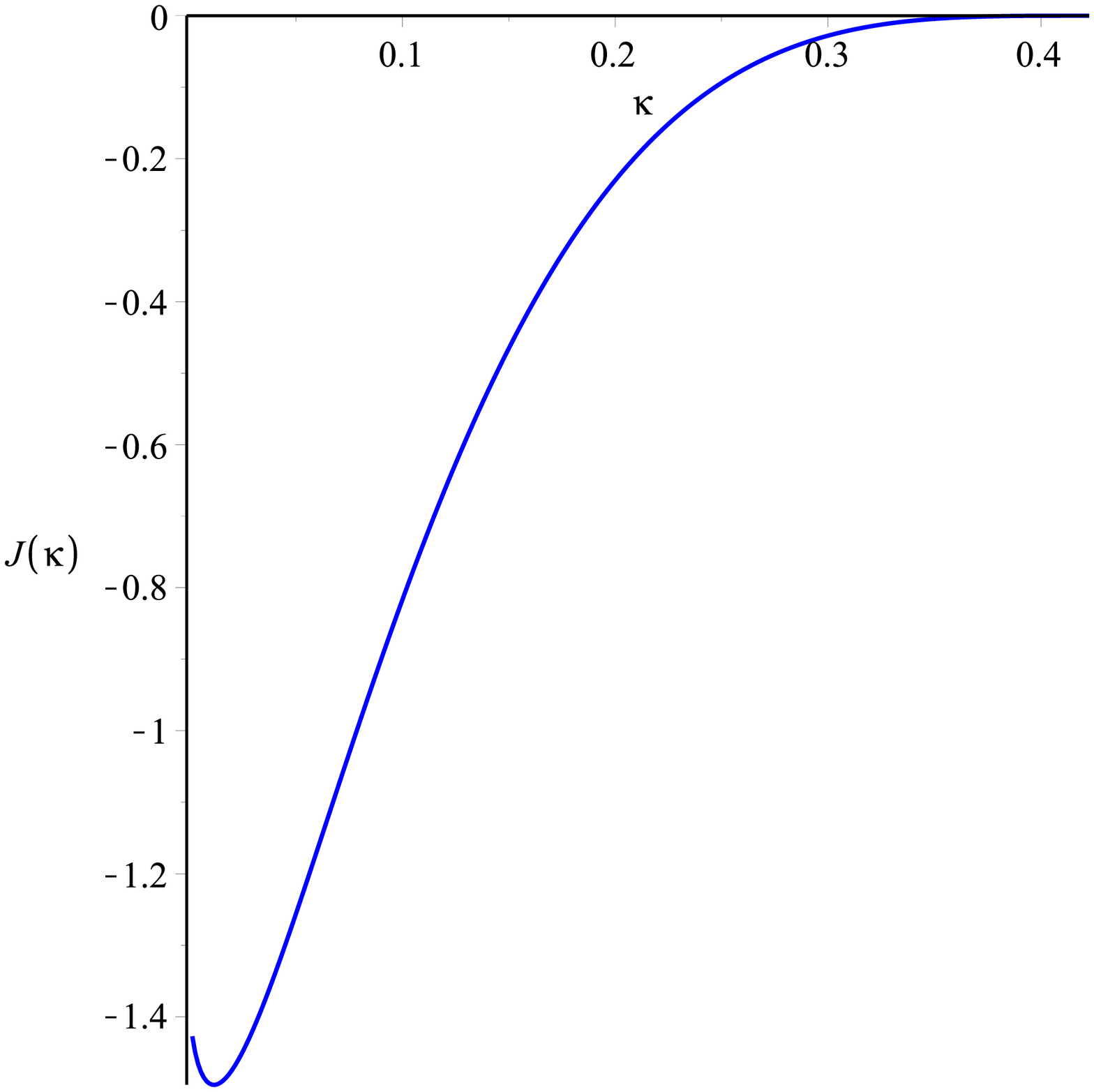}\\
\footnotesize{ (a) (${\kappa}$,\ $I(\kappa)$)} & \footnotesize{ (b) (${\kappa}$,\ $J(\kappa)$)}
\end{tabular}
\end{center}
\begin{center}
\footnotesize {{Fig. 7} \ \  The algebraic curves of $I({\kappa})$ and $J({\kappa})$ with respect to ${\kappa}$ for $0<{\kappa}<1-\frac{\sqrt3}{3}$.}
\end{center}
We can see that $I({\kappa})>0$ and $J({\kappa})<0$ for $1+\frac{\sqrt3}{3}<{\kappa}<2$ in Fig. 6, $I({\kappa})>0$ and $J({\kappa})<0$ for $0<{\kappa}<1-\frac{\sqrt3}{3}$ in Fig. 7. It is shown that a good agreement with the proof of \textbf{Lemma} \ref{lemma1}.

Hence, we substitute (\ref{I2}) and (\ref{J2}) into (\ref{Melnikov}), and the expression of the homoclinic Melnikov function of system (\ref{ODE6}) for $1+\frac{\sqrt3}{3}<{\kappa}<2$ can be rewritten as
\begin{equation}\label{Hom-Melnikov}
\begin{array}{rcl}
&&M_{\textmd{hom}}(c,{\kappa})\Big|_{\kappa\in(1+\frac{\sqrt3}{3},2)}
\\&=&
\frac{1}{c^2}\frac{\sqrt3}{3}\big[6\textmd{ln}2{{\kappa}}^3+6{{\kappa}}^3\textmd{ln}(\sqrt{6{{\kappa}}^2-12{\kappa}+4}+2{\kappa}-2)-18\textmd{ln}2{{\kappa}}^2\\
&&-18{{\kappa}}^2\textmd{ln}(\sqrt{6{{\kappa}}^2-12{\kappa}+4}+2{\kappa}-2)+12\textmd{ln}2{\kappa}\\
&&+12{\kappa}\textmd{ln}(\sqrt{6{{\kappa}}^2-12{\kappa}+4}+2{\kappa}-2)+2\sqrt{6{{\kappa}}^2-12{\kappa}+4}~ \big]\\
&&-\frac{1}{c^2}\sqrt3 {\kappa}\big[{{\kappa}}^2\textmd{ln}(2{\kappa}-{{\kappa}}^2)+3\textmd{ln}2{{\kappa}}^2-3{\kappa}\textmd{ln}(2{\kappa}-{{\kappa}}^2)-9\textmd{ln}2{\kappa}\\
&&+2\textmd{ln}(2{\kappa}-{{\kappa}}^2)+6\textmd{ln}2\big]-\frac{\sqrt3}{5} \big[27{{\kappa}}^4 \sqrt{6{{\kappa}}^2-12{\kappa}+4}\\
&&+15 {{\kappa}}^3\textmd{ln}(2{\kappa}-{{\kappa}}^2)-198{{\kappa}}^3\sqrt{6{{\kappa}}^2-12{\kappa}+4}\\
&&-30{{\kappa}}^3\textmd{ln}(\sqrt{6{{\kappa}}^2-12{\kappa}+4}+2{\kappa}-2)+15\textmd{ln}2{{\kappa}}^3-45{{\kappa}}^2\textmd{ln}(2{\kappa}-{{\kappa}}^2) \\
&&+90{{\kappa}}^2\textmd{ln}(\sqrt{6{{\kappa}}^2-12{\kappa}+4}+2{\kappa}-2)+30\sqrt2{\kappa}\sqrt{(3{{\kappa}}^2-6{\kappa}+2)^3}\\
&&+294{{\kappa}}^2\sqrt{6{{\kappa}}^2-12{\kappa}+4}-45\textmd{ln}2{{\kappa}}^2+30{\kappa}\textmd{ln}(2{\kappa}-{{\kappa}}^2)\\
&&-60{\kappa}\textmd{ln}(\sqrt{6{{\kappa}}^2-12{\kappa}+4}+2{\kappa}-2)+10\sqrt2\sqrt{(3{{\kappa}}^2-6{\kappa}+2)^3}\\
&&-72{\kappa}\sqrt{6{{\kappa}}^2-12{\kappa}+4}+30\textmd{ln}2{\kappa}-18\sqrt{6{{\kappa}}^2-12{\kappa}+4}~ \big].
\end{array}
\end{equation}
{It is easy to verify that for $\kappa\in(1+\frac{\sqrt3}{3},2)$, there is $c=c_1(\kappa)$ such that $M_{\textmd{hom}}(c,{\kappa})\Big|_{\kappa\in(1+\frac{\sqrt3}{3},2)}=0$.}

Similarly, the expression of the homoclinic Melnikov function of system (\ref{ODE6}) for $0<{\kappa}<1-\frac{\sqrt3}{3}$ is
\begin{equation}\label{Hom-Melnikov1}
\begin{array}{rcl}
&&
M_{\textmd{hom}}(c,{\kappa})\Big|_{\kappa\in(0,1-\frac{\sqrt3}{3})}
\\
&=&\frac{1}{c^2}\frac{\sqrt3}{3}\big[6\textmd{ln}2{{\kappa}}^3+6{{\kappa}}^3\textmd{ln}(\sqrt{6{{\kappa}}^2-12{\kappa}+4}+2{\kappa}-2)-18\textmd{ln}2{{\kappa}}^2\\
&&-18{{\kappa}}^2\textmd{ln}(\sqrt{6{{\kappa}}^2-12{\kappa}+4}+2{\kappa}-2)+12\textmd{ln}2{\kappa}\\
&&+12{\kappa}\textmd{ln}(\sqrt{6{{\kappa}}^2-12{\kappa}+4}+2{\kappa}-2)+2\sqrt{6{{\kappa}}^2-12{\kappa}+4}~ \big]\\
&&-\frac{1}{c^2}\sqrt3 {\kappa}\big[2{{\kappa}}^2\textmd{ln}(-\sqrt{2{\kappa}-{{\kappa}}^2})+3\textmd{ln}2{{\kappa}}^2-6{\kappa}\textmd{ln}(-\sqrt{2{\kappa}-{{\kappa}}^2})-9\textmd{ln}2{\kappa}\\
&&+4\textmd{ln}(-\sqrt{2{\kappa}-{{\kappa}}^2})+6\textmd{ln}2\big]-\frac{\sqrt3}{5} \big[27{{\kappa}}^4 \sqrt{6{{\kappa}}^2-12{\kappa}+4}\\
&&+30 {{\kappa}}^3\textmd{ln}(-\sqrt{2{\kappa}-{{\kappa}}^2})+60{\kappa}\textmd{ln}(-\sqrt{2{\kappa}-{{\kappa}}^2})+15\textmd{ln}2{{\kappa}}^3\\
&&-30{{\kappa}}^3\textmd{ln}(\sqrt{6{{\kappa}}^2-12{\kappa}+4}+2{\kappa}-2)-90{{\kappa}}^2\textmd{ln}(-\sqrt{2{\kappa}-{{\kappa}}^2}) \\
&&+90{{\kappa}}^2\textmd{ln}(\sqrt{6{{\kappa}}^2-12{\kappa}+4}+2{\kappa}-2)+30\sqrt2{\kappa}\sqrt{(3{{\kappa}}^2-6{\kappa}+2)^3}\\
&&+294{{\kappa}}^2\sqrt{6{{\kappa}}^2-12{\kappa}+4}-45\textmd{ln}2{{\kappa}}^2-198{{\kappa}}^3\sqrt{6{{\kappa}}^2-12{\kappa}+4}\\
&&-60{\kappa}\textmd{ln}(\sqrt{6{{\kappa}}^2-12{\kappa}+4}+2{\kappa}-2)+10\sqrt2\sqrt{(3{{\kappa}}^2-6{\kappa}+2)^3}\\
&&-72{\kappa}\sqrt{6{{\kappa}}^2-12{\kappa}+4}+30\textmd{ln}2{\kappa}-18\sqrt{6{{\kappa}}^2-12{\kappa}+4}\ \big].
\end{array}
\end{equation}
{It is easy to verify that for $\kappa\in(0,1-\frac{\sqrt3}{3})$, there is $c=c_2(\kappa)$ such that $M_{\textmd{hom}}(c,{\kappa})\Big|_{\kappa\in(0,1-\frac{\sqrt3}{3})}=0$.}

Secondly, the Melnikov method is also employed to detect the existence of kink (anti-kink) wave solutions under small perturbation. The heteroclinic Melnikov function of system (\ref{ODE6}) is defined as follows:
\begin{equation}\label{Melnikov-het}
\begin{array}{rcl}
M_{\textmd{het}}(c,\kappa)&=&\displaystyle\int_{\Upsilon(\kappa)}(\frac{9}{2}z^2-9z+3+\frac{1}{c^2})y^2d\zeta=\displaystyle\int_{-\infty}^{+\infty}(\frac{9}{2}z^2-9z+3+\frac{1}{c^2})y^2d\zeta\\
&=&\displaystyle\int_{0}^{2}(\frac{9}{2}z^2-9z+3+\frac{1}{c^2})ydz\\
&=&\frac{\sqrt3}{2}\displaystyle\int_{0}^{2}(\frac{9}{2}z^2-9z+3+\frac{1}{c^2})z(2-z)dz=\displaystyle-\frac{2\sqrt3(3c^2-5)}{15c^2}.
\end{array}
\end{equation}
\begin{mylem}\label{lemma2}
For $\kappa=0$ or $\kappa=2$, $M_{het}(c,\kappa)=0$ has a positive root $c^\ast=\frac{\sqrt{15}}{3}$. Moreover, $\frac{\partial M_{\textmd{het}}}{\partial c}\mid_{c=c^\ast}\neq0$.
\end{mylem}
\emph{Proof. } Clearly, it directly follows form (\ref{Melnikov-het}) that $M_{\textmd{het}}(c,\kappa)=0$ has a positive root $c^\ast=\frac{\sqrt{15}}{3}$. Meanwhile,
$$\frac{\partial M_{\textmd{het}}}{\partial c}= -\frac{4\sqrt{3}}{5c}+\frac{4\sqrt{3}(3c^2-5)}{15c^3},$$
and  $$\frac{\partial M_{\textmd{het}}}{\partial c}\mid_{c=c^\ast}=-\frac{12\sqrt{5}}{25}\neq0.$$

\noindent{\bf Proof of the Theorem \ref{theorem1}:} By Lemma \ref{lemma1} and Lemma \ref{lemma2}, Theorem \ref{theorem1} follows immediately.

\section{Numerical analysis}
Numerical simulations are employed to confirm the theoretical results derived in previous sections. Here, maple software $18.0$ is used.

Firstly, we simulate the existence of solitary wave solution of perturbed $(1+1)$-dimensional DLWE (\ref{wuzhang1}). Let ${\kappa}=1.7$ and $\epsilon=0.01$. By  (\ref{I2}) and (\ref{J2}), we obtain $I(1.7)\approx0.0648882124$ and $J(1.7)\approx-0.02790800120$ such that $c(1.7)=\sqrt{-\frac{I(1.7)}{J(1.7)}}\approx1.524819859$. Then, according to \textbf{Theorem} \ref{theorem1}, we take $c=c(1.7)$ and let the initial value be $(z(0),y(0))=(z_+,0)=(1.309950494,0)$ (red point) which the homoclinic orbit would pass through. The phase portraits $(z,y)$ and time history curves $(\zeta,z)$ of system (\ref{ODE6}) are plotted in Fig. 8.
\begin{center}
\begin{tabular}{ccc}
\epsfxsize=5.5cm \epsfysize=5.5cm \epsffile{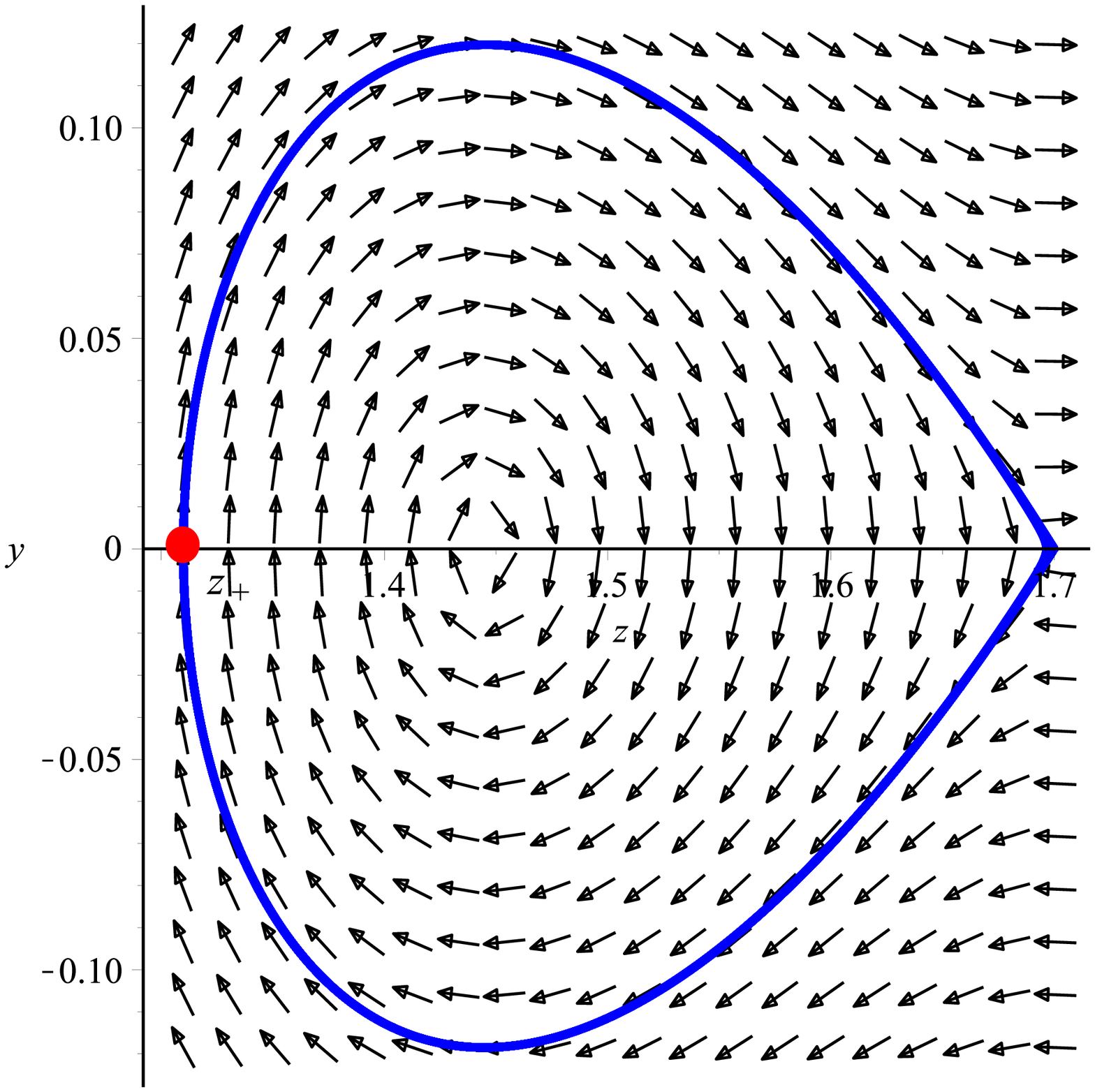}&
\epsfxsize=5.5cm \epsfysize=5.5cm \epsffile{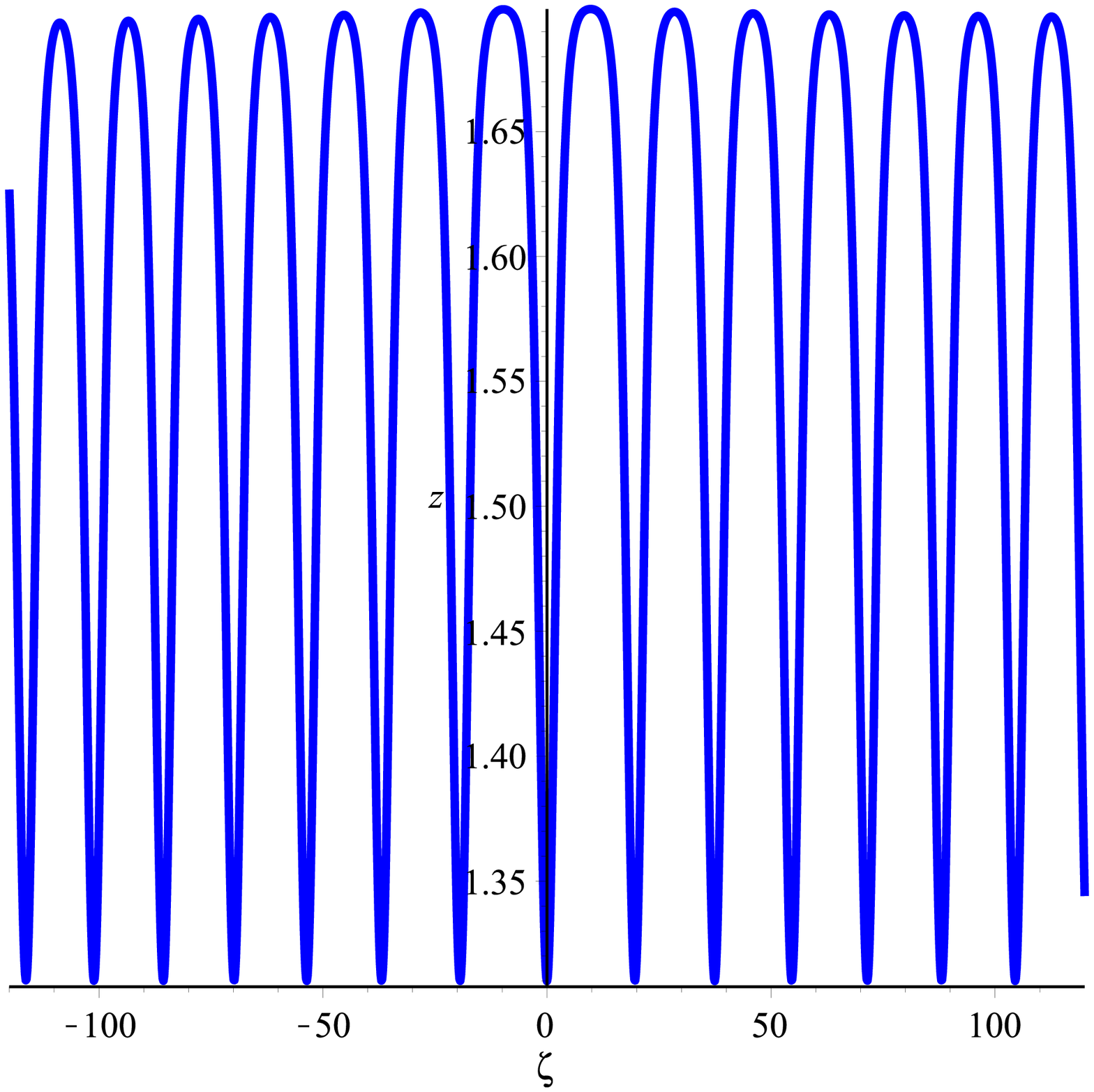}\\
\footnotesize{ (a) $c=c(1.7)$ } & \footnotesize{ (b) $(\zeta,z)$ }
\end{tabular}
\end{center}
\begin{center}
\footnotesize {{Fig. 8} \ (a) Phase portraits with $c=c(1.7)$, (b) Time history curves of $(\zeta,z)$ for system (\ref{ODE6}) for $\epsilon=0.01$, ${\kappa}=1.7$ and initial value $(z(0),y(0))=(z_+,0)=(1.309950494,0)$.}
\end{center}
One can see that the homoclinic orbit and dark solitary wave solution of system (\ref{ODE6}) still persist after a small singular perturbation.

Secondly, set ${\kappa}=0.3$ and $\epsilon=0.01$. By (\ref{I22}) and (\ref{J22}), we have $I(0.3)\approx0.0648882229$ and $J(0.3)\approx-0.02790784012$ such that $c(0.3)=\sqrt{-\frac{I(0.3)}{J(0.3)}}\approx1.524824377$. Taking $c=c(0.3)$ and initial value to be $(z(0),y(0))=(z_-,0)=(0.690049506,0)$ (red point) which the homoclinic orbit would pass through. The phase portraits $(z,y)$ and time history curves $(\zeta,z)$ of system (\ref{ODE6}) are plotted in Fig. 9.
\begin{center}
\begin{tabular}{ccc}
\epsfxsize=5.5cm \epsfysize=5.5cm \epsffile{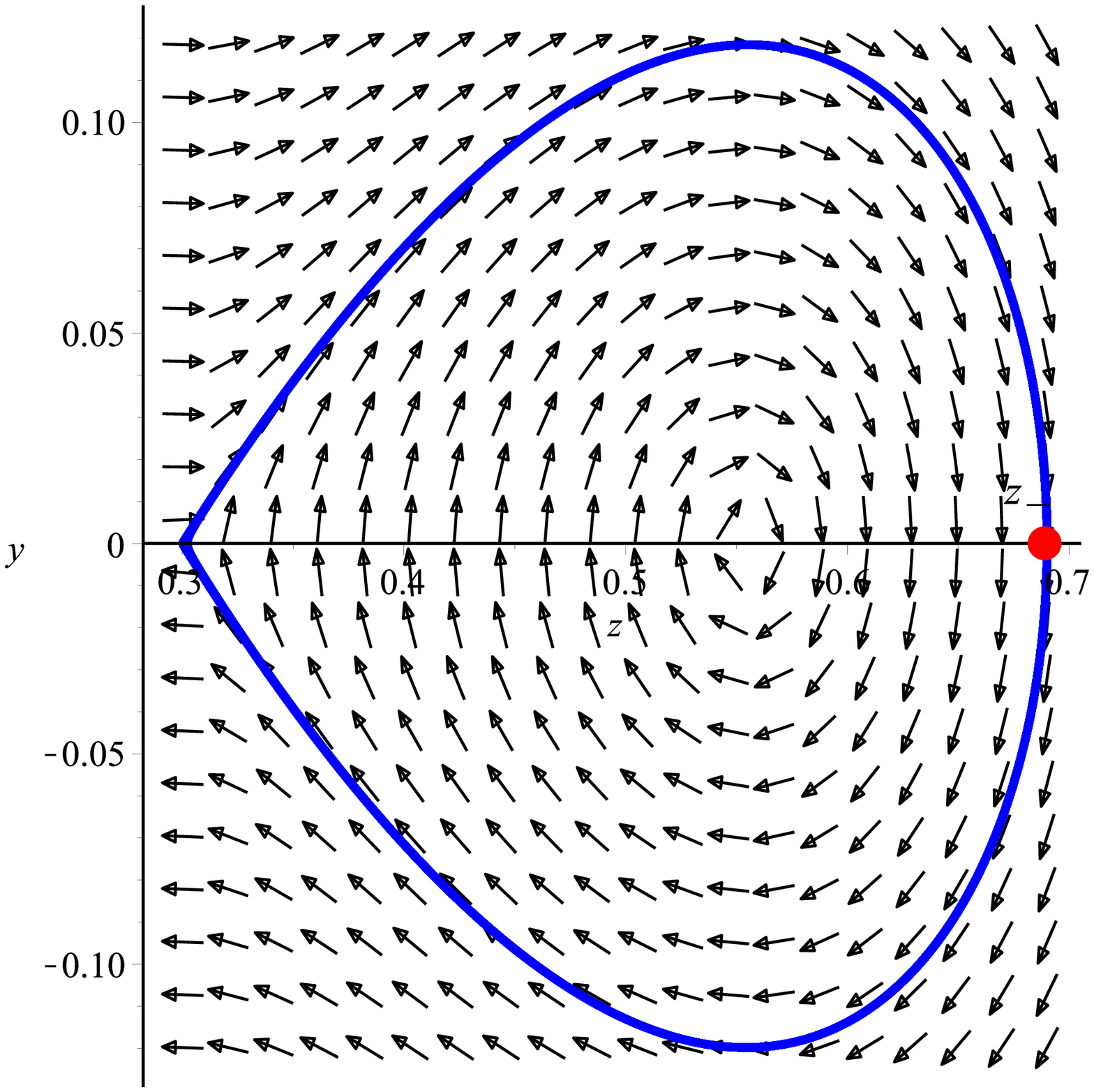}&
\epsfxsize=5.5cm \epsfysize=5.5cm \epsffile{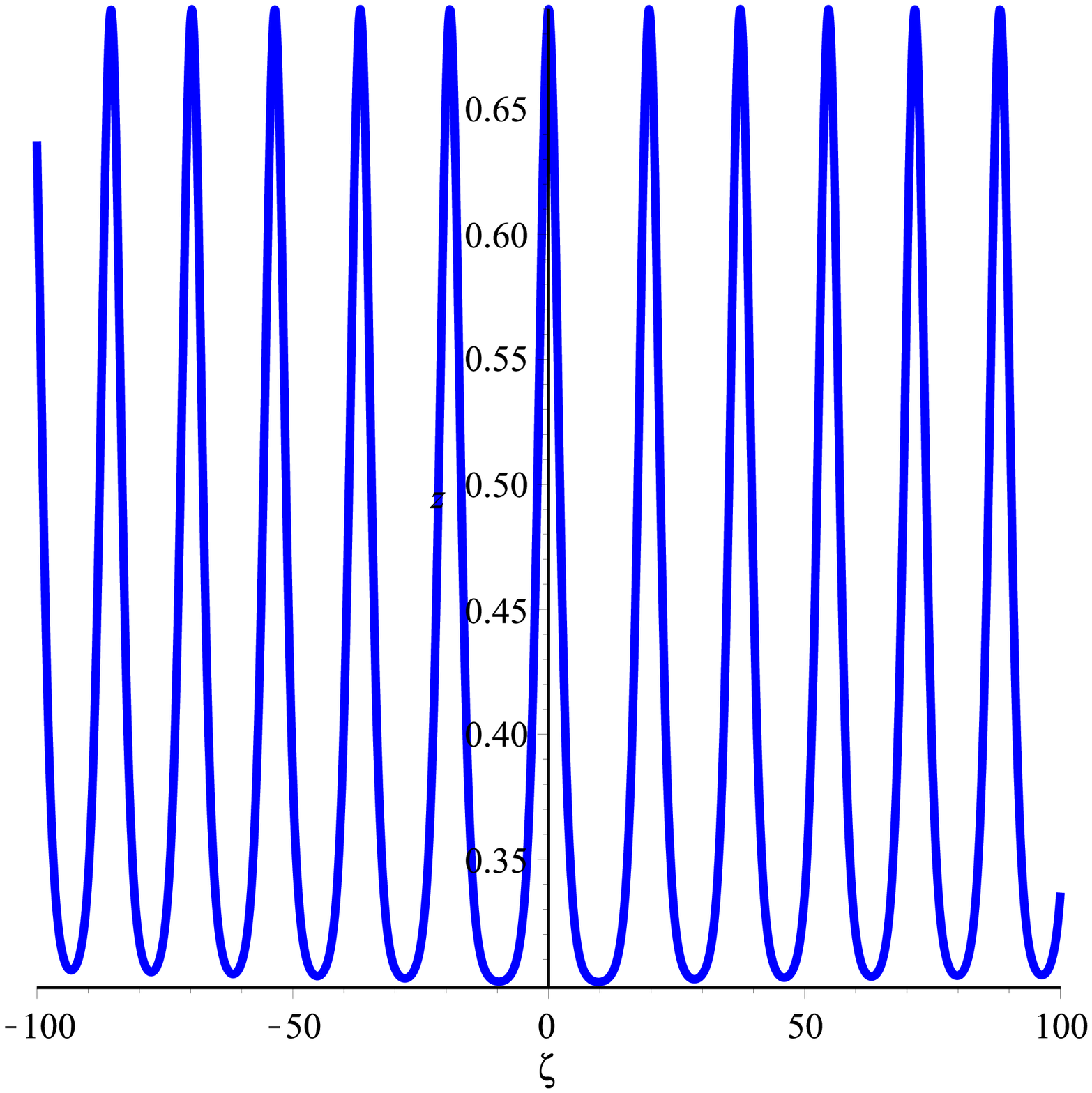}\\
\footnotesize{ (a) $c=c(0.3)$ } & \footnotesize{ (b) $(\zeta,z)$ }
\end{tabular}
\end{center}
\begin{center}
\footnotesize {{Fig. 9} \  (a) Phase portraits with $c=c(0.3)$, (b) Time history curves of $(\zeta,z)$ for system (\ref{ODE6}) for $\epsilon=0.01$, ${\kappa}=0.3$ and initial value $(z(0),y(0))=(z_-,0)=(0.690049506,0)$.}
\end{center}
We see that the homoclinic orbit and bright solitary wave solution of system (\ref{ODE6}) still exist after a small singular perturbation which verifies the \textbf{Theorem} \ref{theorem1} $(i)$.

Next, let us to see the existence of kink and anti-kink wave solution of system (\ref{ODE6}). Giving $\kappa=0$ or $\kappa=2$, $\epsilon=0.01$ and initial value $(z(0),y(0))=(1,\pm\frac{\sqrt3}{2})$ (red point). We draw the phase portraits of system (\ref{ODE6}) after taking $c=\frac{\sqrt{15}}{3}$ in Fig. 10. Form the Fig. 10 (a), it is shown that the heteroclinic orbits still persist. Correspondingly, a pair of kink and anti-kink wave solutions exists (see Fig. 10 (b) and (c)), which is consistent with \textbf{Theorem} \ref{theorem1} $(ii)$.
\begin{center}
\begin{tabular}{ccc}
\epsfxsize=4.5cm \epsfysize=4.5cm \epsffile{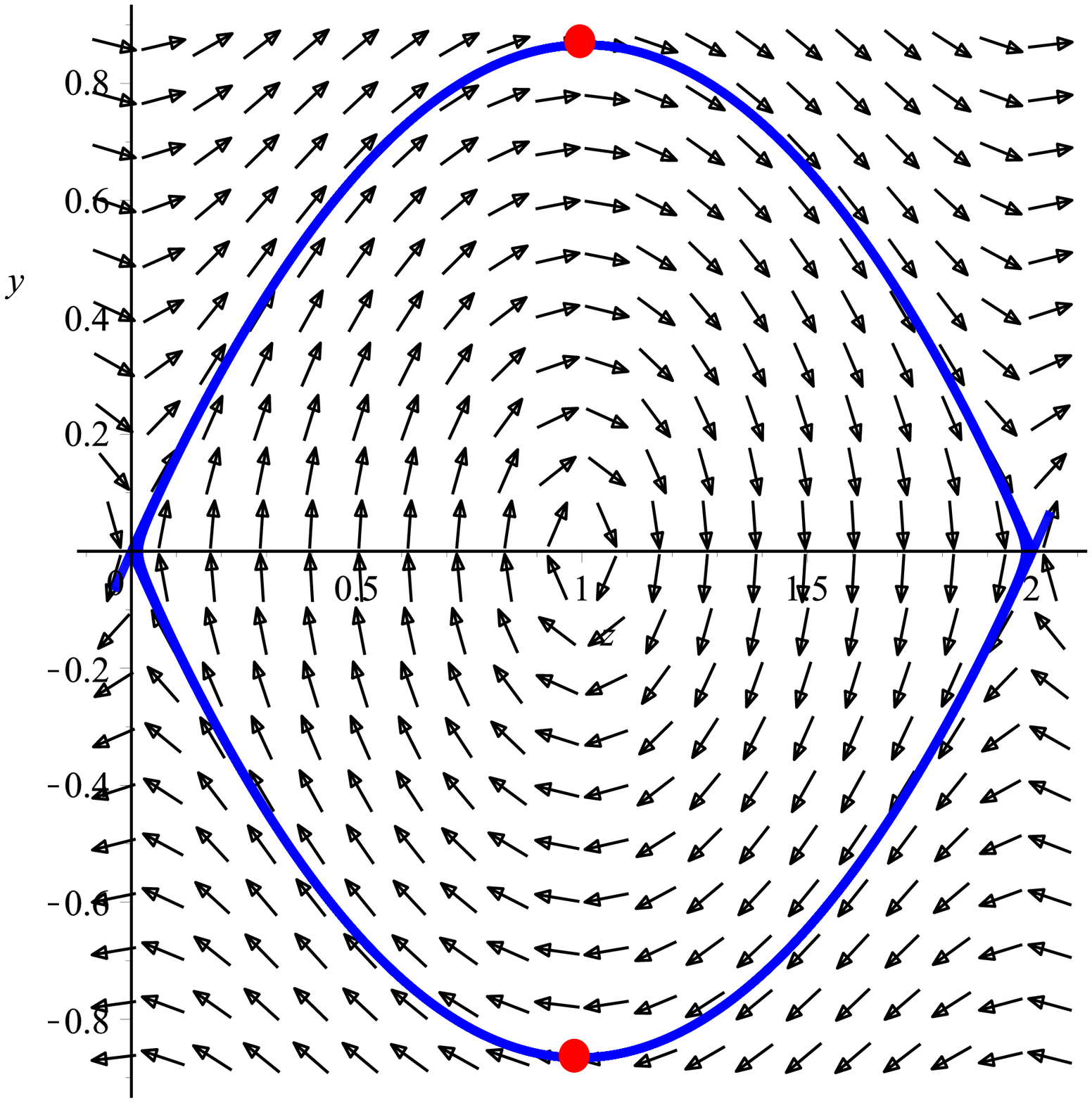}&
\epsfxsize=4.5cm \epsfysize=4.5cm \epsffile{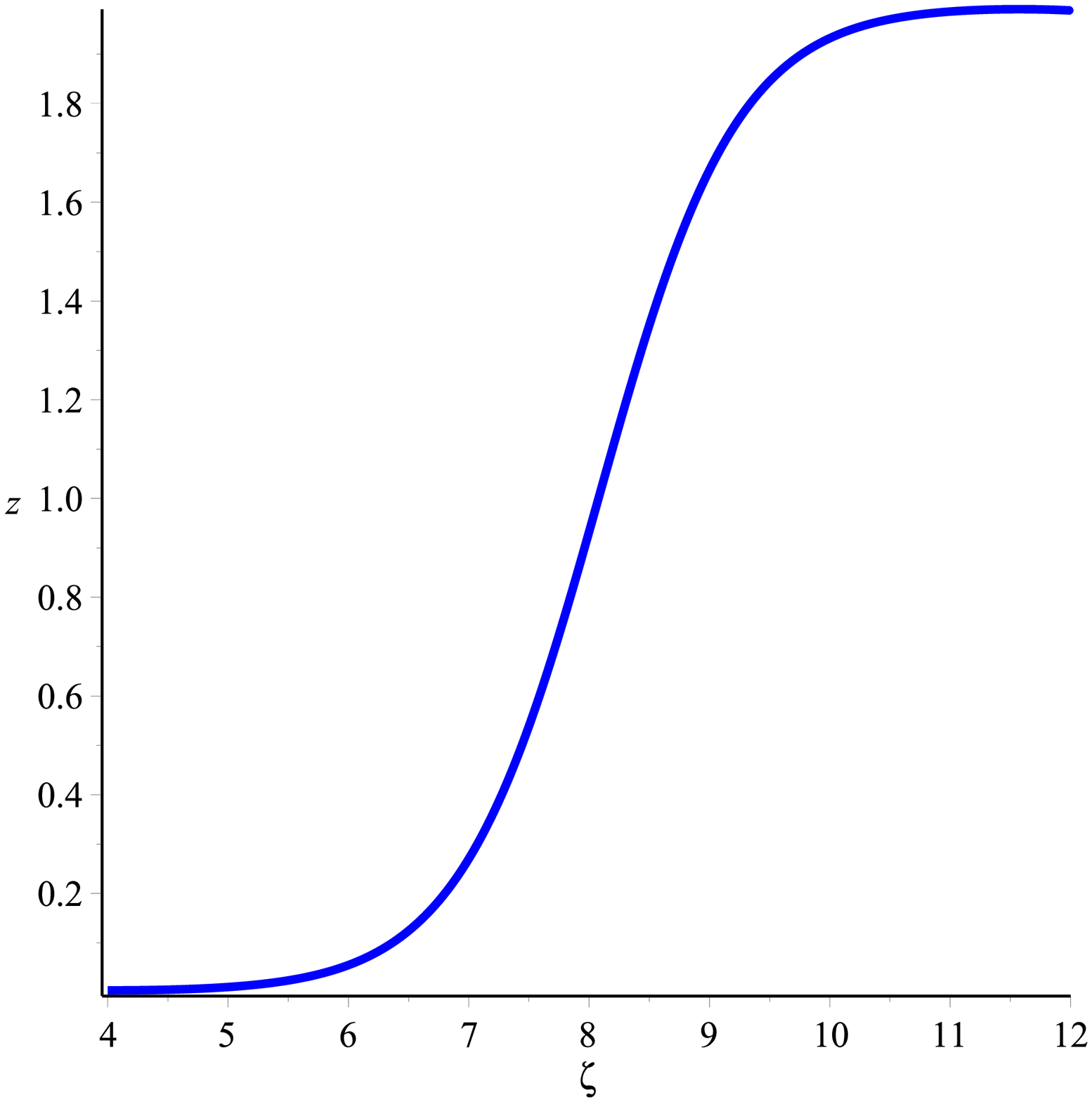}&
\epsfxsize=4.5cm \epsfysize=4.5cm \epsffile{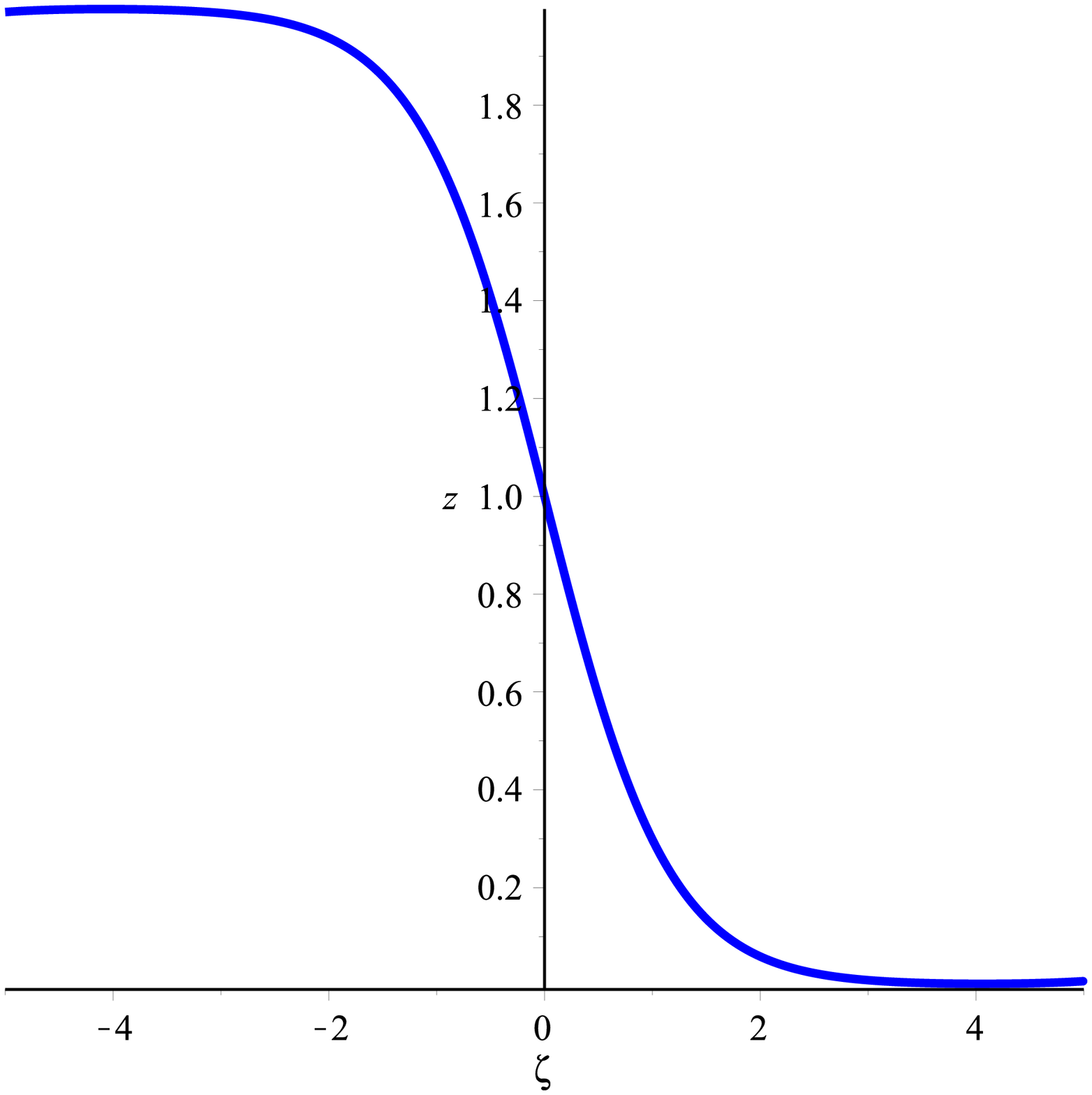}\\
\footnotesize{ (a) $c=\frac{\sqrt{15}}{3}$ } & \footnotesize{ (b) Kink wave} & \footnotesize{ (c) Anti-kink wave}
\end{tabular}
\end{center}
\begin{center}
\footnotesize {{Fig. 10} \ Phase portraits with $c=\frac{\sqrt{15}}{3}$ for $\epsilon=0.01$, initial value $(z(0),y(0))=(1,\pm\frac{\sqrt3}{2})$.}
\end{center}

%\section{Conclusion}
%{This paper mainly studies the existence of solitary and kink wave solutions for a perturbed $(1+1)$-dimensional DLWE by using GSP approach, Melnikov method and bifurcation theorem. We obtain the exact solitary and kink (anti-kink) wave solutions of unperturbed $(1+1)$-dimensional DLWE, then we prove that the solitary wave solutions exist with a suitable $c=c({\kappa})$ and kink (anti-kink) wave solutions exist with a unique wave speed $c^*$. Whereafter, with the help of Maple, numerical simulations illustrate the previous theoretical results.}

\section{Fund Acknowledgement}
This work was jointly supported by the Zhejiang Provincial Natural Science Foundation of China (No.LZ23A010001), National Natural Science Foundation of China under Grant (No. 11671176, 11931016), Natural Science Foundation of Fujian Province under Grant (No. 2021J011148), Fujian Province Young Middle-Aged teachers education scientific research project (No. JAT210454) and Teacher and Student Scientific Team Fund of Wuyi University (Grant No. 2020-SSTD-003).

\section{ Conflict of Interest}
 The authors declare that they have no conflict of interest.

\section{Data Availability Statement}

No data was used for the research in this article.

\section{ Author contributions}
 All the authors have same contributions to the paper.

%\end{multicols}
\end{document}